%% file: paper.tex
\documentclass{siamltex}
\usepackage[numbers]{natbib}
\usepackage{amsmath}
\usepackage{amscd}
\usepackage{amssymb}
\usepackage{graphicx}
\usepackage{url}
\usepackage{color}
\usepackage{listings}
\usepackage{booktabs}
\usepackage{textcomp}
\usepackage{mathtools}
\usepackage{hyperref}
\usepackage[utf8]{inputenc}
\usepackage{tikz}
\usepackage{pgfplots}
\usepackage{environ}
\usepackage{algorithm}
\usepackage{enumerate}
\usepackage[noend]{algpseudocode}
\usepackage{subfigure}

\usepackage{color}
\usepackage{todonotes}

\newcommand{\sech}{\mathrm{sech}}

\newcommand{\ov}{\overline}

\newcommand{\non}{\nonumber}
\newcommand{\eps}{\epsilon}
\newcommand{\beqa}{\begin{eqnarray}}
\newcommand{\eeqa}{\end{eqnarray}}
\newcommand{\beqas}{\begin{eqnarray*}}
\newcommand{\eeqas}{\end{eqnarray*}}
\newcommand{\ba}{\begin{align}}
\newcommand{\ea}{\end{align}}
\newcommand{\beq}{\begin{equation}}
\newcommand{\eeq}{\end{equation}}

\newcommand{\ra}{\rightarrow}

\newcommand{\Z}{{\mathbb Z}}

\newcommand{\pdhfrac}[2]{\mathchoice{\frac{#1}{#2}}{#1/#2}{#1/#2}{#1/#2}}
\newcommand{\fdd}[2]{\pdhfrac{\mathrm{d}#1}{\mathrm{d}#2}}
\newcommand{\sdd}[2]{\pdhfrac{\mathrm{d}^2#1}{\mathrm{d}#2^2}}
\newcommand{\pd}[2]{\pdhfrac{{\partial}#1}{{\partial}#2}}
\newcommand{\spd}[2]{\pdhfrac{\partial^2#1}{{\partial}#2^2}}
\newcommand{\mpd}[3]{\pdhfrac{\partial^2#1}{{\partial}#2{\partial}#3}}
\newcommand{\nummax}{M}
\newcommand{\maxmax}{M_{\rm max}}
\renewcommand{\d}[1]{\mathrm{d}#1}
\graphicspath{{./figures/}}

\definecolor{blue}{rgb}{0.2980392156862745, 0.4470588235294118, 0.6901960784313725}
\definecolor{green}{rgb}{0.3333333333333333, 0.6588235294117647, 0.40784313725490196}
\definecolor{red}{rgb}{0.7686274509803922, 0.3058823529411765, 0.3215686274509804}

\definecolor{darkblue}{rgb}{0.00,0.00,0.55}
\definecolor{black}{rgb}{0.00,0.00,0.00}
\hypersetup{
    linkcolor = darkblue,
    anchorcolor = darkblue,
    citecolor = darkblue,
    filecolor = darkblue,
    urlcolor = darkblue,
    colorlinks  = true,
}

\title{Analysis of Carrier's problem}

\author{
  S.~J.~Chapman\thanks{Mathematical Institute, University of Oxford, Oxford, UK (\texttt{chapman@maths.ox.ac.uk}).}
  \and
  P.~E.~Farrell\thanks{Mathematical Institute, University of Oxford, Oxford, UK.
    Center for Biomedical Computing, Simula Research Laboratory, Oslo, Norway
    (\texttt{patrick.farrell@maths.ox.ac.uk}).
    This research is supported by EPSRC grants EP/K030930/1 and EP/M019721/1, by a Center of Excellence grant from the Research
    Council of Norway to the Center for Biomedical Computing at Simula
    Research Laboratory (project number 179578), and by the generous support of Sir Michael Moritz and
    Harriet Heyman. The authors would like to thank \'A.~Birkisson and L.~N.~Trefethen for useful discussions.}
  }

\begin{document}
\maketitle

\begin{abstract}
A computational and asymptotic analysis of the solutions of Carrier's problem is
presented. The computations reveal a striking and beautiful bifurcation diagram,
with an infinite sequence of alternating pitchfork and fold bifurcations as
the bifurcation parameter tends to zero. The method of
Kuzmak is then applied to construct 
asymptotic solutions to the problem. This asymptotic approach explains
the bifurcation 
structure identified numerically, and its predictions of
the bifurcation points are in excellent agreement with the
numerical results. The analysis yields a 
novel and complete taxonomy of the solutions to the problem, and demonstrates
that a claim of Bender \& Orszag \cite{bender1999} is incorrect. 
\end{abstract}

\begin{keywords}
Multiple scales, Kuzmak's method, bifurcation, asymptotic analysis.
\end{keywords}

\begin{AMS}
 34E13, 41A60, 34E05
\end{AMS}

\section{Introduction}
\label{intro}
In 1970, G.~F.~Carrier \cite[eq.~(3.5)]{carrier1970} introduced
the following singular perturbation problem
\beq
 \eps^2 y'' + 2(1-x^2) y + y^2 = 1, \qquad y(-1) = y(1) = 0,\label{maineqn}
\eeq
where $0< \eps \ll 1$, and a prime represents $\fdd{}{x}$.
This remarkably beautiful and complex problem is discussed in more
detail in the textbooks of Carrier \& Pearson
\cite[p.~197]{carrier1985} and Bender \& Orszag
\cite[p.~464]{bender1999}. We briefly review their discussion.

\begin{figure}
\centering
\subfigure[]{
{\includegraphics[width=0.45\textwidth]{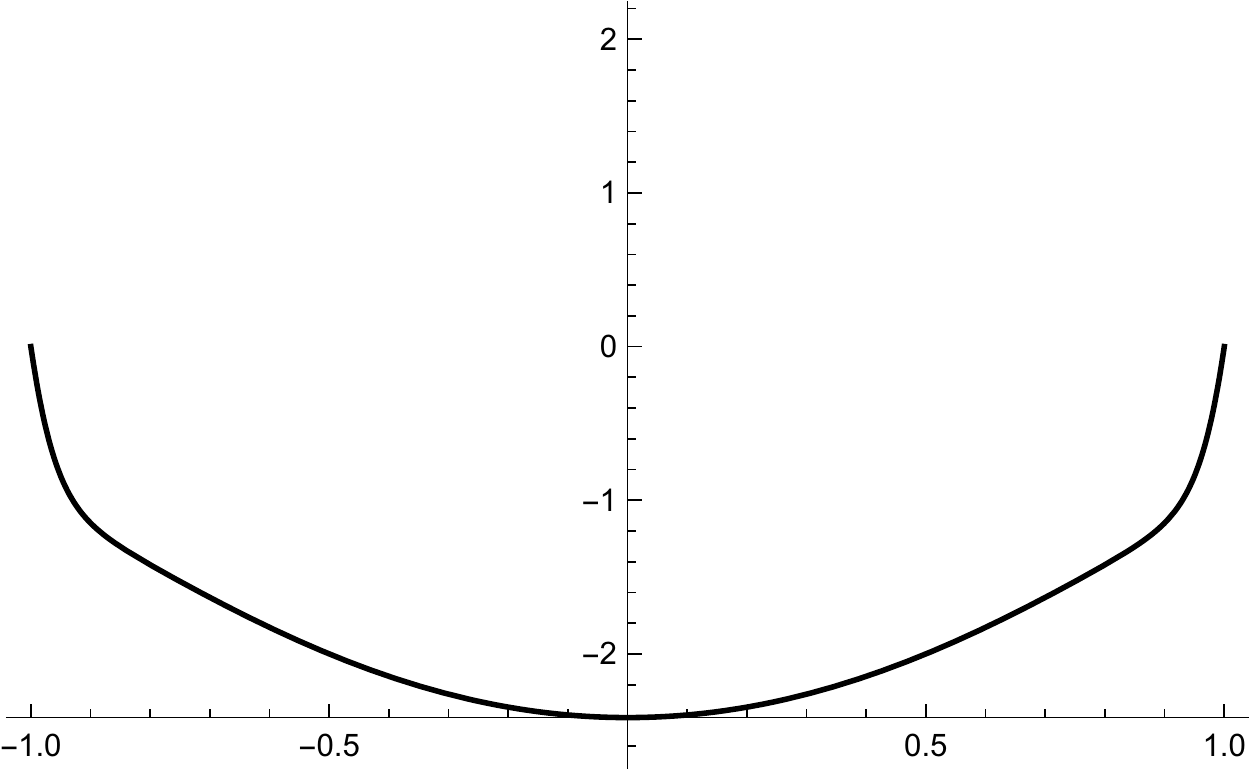}}
}
\
\subfigure[]{
{\includegraphics[width=0.45\textwidth]{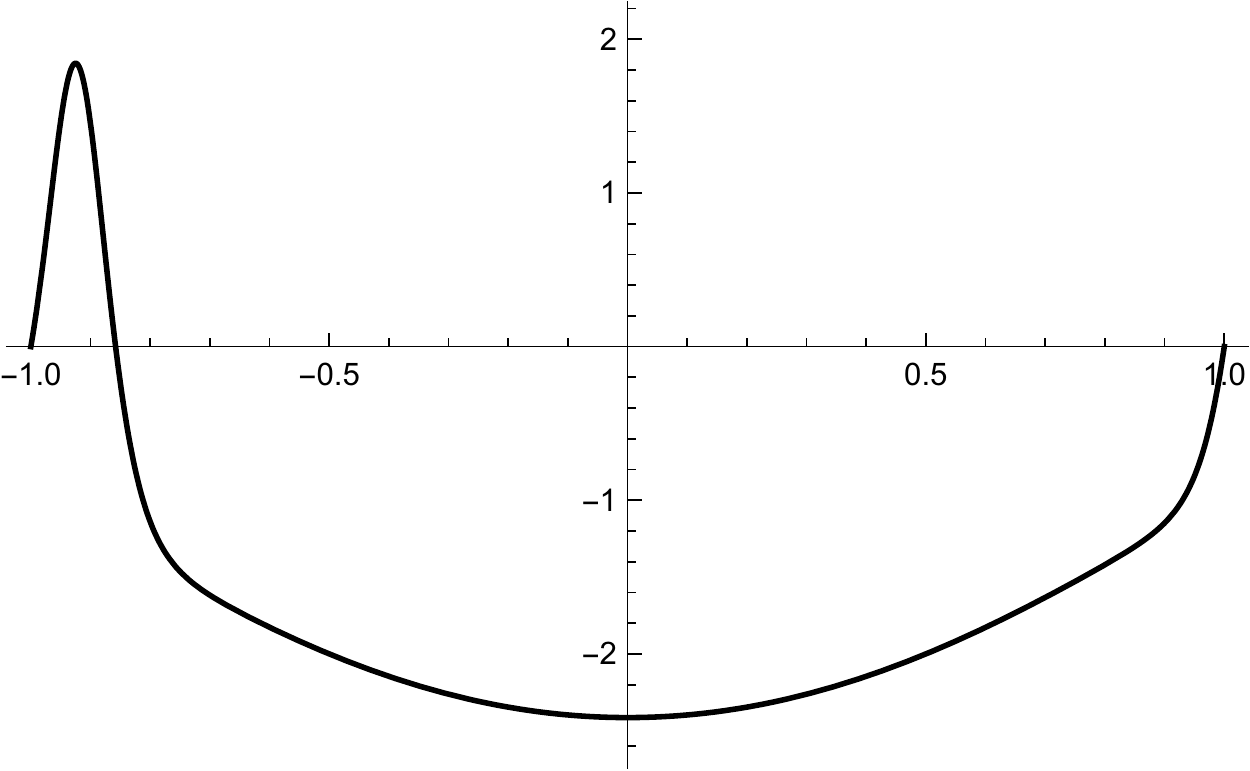}}
}
\\
\subfigure[]{
{\includegraphics[width=0.45\textwidth]{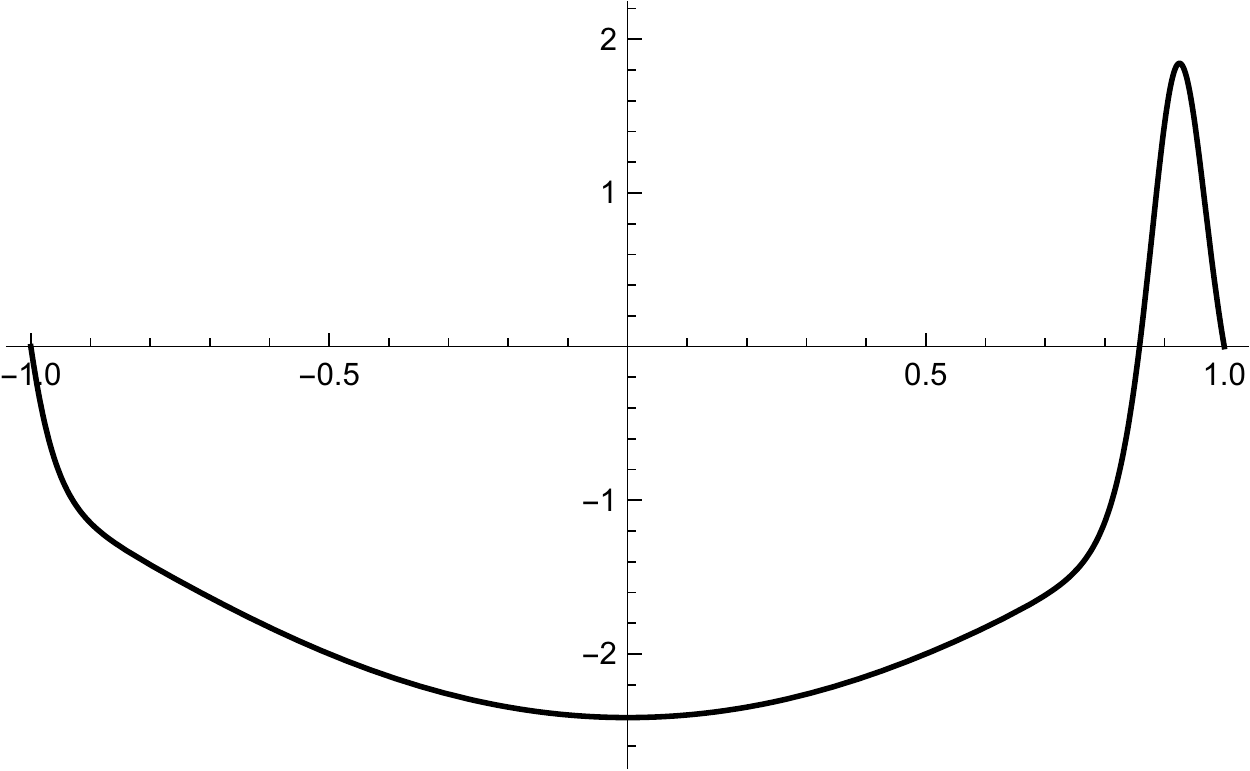}}
}
\
\subfigure[]{
{\includegraphics[width=0.45\textwidth]{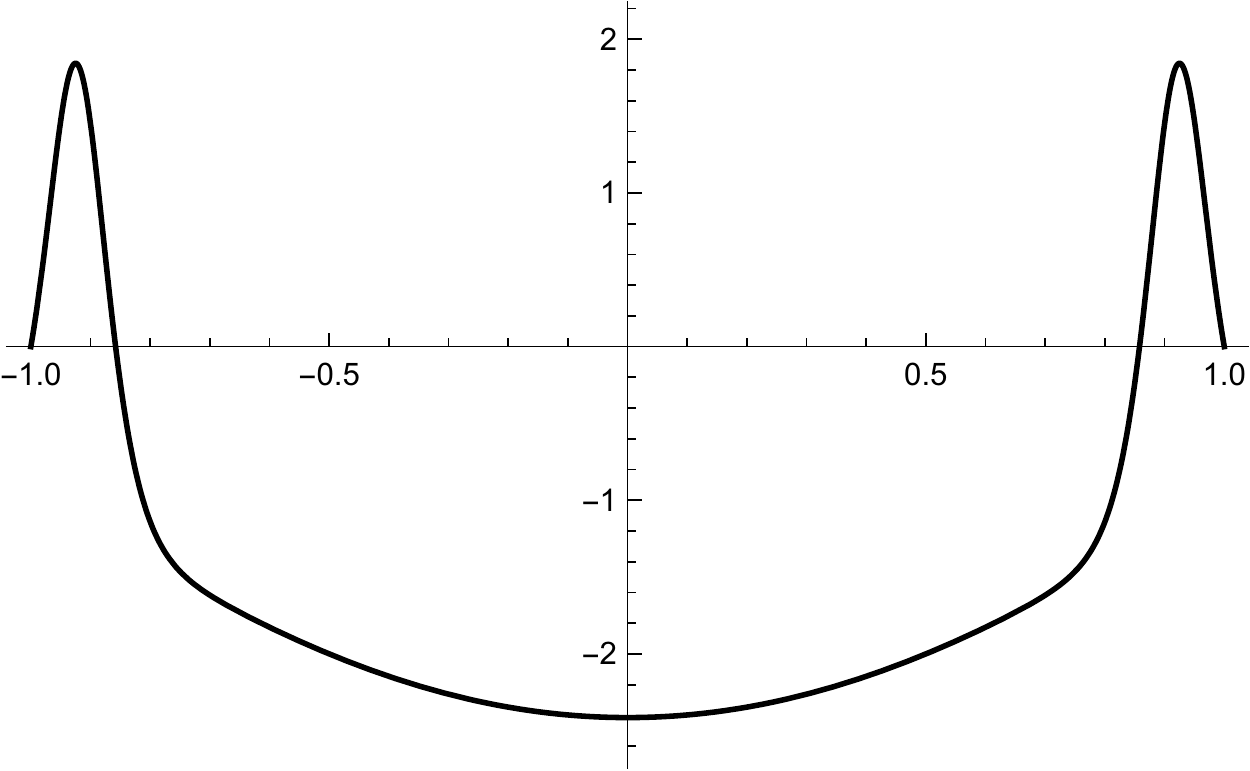}}
}
\caption{The four asymptotic solutions with $\eps^2 = 0.00223$, generated
  using the outer approximation (\ref{outersol}) and the boundary layer
  approximations (\ref{inner}). }
\label{fig1}
\end{figure}

Since \eqref{maineqn} is singularly perturbed, we expect the solution
to comprise an outer solution valid for $1-|x| \gg \eps$ combined with
possible boundary layers near $x = \pm 1$.
Na\"ively setting $\eps = 0$ gives the leading-order outer solution as
\beq
 y_{\mathrm{out}} = x^2-1 \pm \sqrt{1+(1-x^2)^2}.\label{outersol}
\eeq
For neither choice of sign does $y_{\mathrm{out}}$ satisfy the
boundary conditions at $x = \pm 1$, so there are indeed boundary layers.
In the boundary layer near $x=-1$ we set $x = -1 + \eps X$,
$y(x) = y_{\mathrm{in}}(X)$ to give
\beq 
\sdd{y_{\mathrm{in}}}{X} + y_{\mathrm{in}}^2 = 1, \qquad
    y_{\mathrm{in}}(0) = 0.\label{innereqn}
\eeq
In order to match with the outer solution, $y_{\mathrm{in}}$ must tend
to $\pm 1$ as $X \ra \infty$. Bender \& Orszag show that there are no
solutions  tending to $1$, so that the minus sign must be chosen in
the outer approximation (\ref{outersol}). On the other hand, there are
{\em two} solutions of (\ref{innereqn}) which tend to $-1$ at infinity, namely
\beq
y_{\mathrm{in}} = -1 + 3\, \sech^2\left( \pm \frac{X}{\sqrt{2}} +
   \tanh^{-1}\left(\sqrt{\frac{2}{3}}\right) \right).\label{inner}
\eeq
Similarly there are two possible boundary layer solutions near $x =
1$.
Thus it seems that a matched asymptotic analysis has produced four
independent solutions of the equation. 
As Bender \& Orszag say, ``it
is a glorious triumph of boundary layer theory that all four solutions
actually exist and are extremely well approximated by the
leading-order uniform approximation'' generated from (\ref{outersol}) and
(\ref{inner}).  These uniform approximations are shown in Figure
\ref{fig1}. 

However, the story does not end there. Bender \& Orszag show that it
is also possible to have a solution with an internal  layer
near $x=0$. Writing $x = \eps X$, $y(x) = y_{\mathrm{in}}(X)$ gives
\beq 
\sdd{y_{\mathrm{in}}}{X} + 2 y_{\mathrm{in}}+ y_{\mathrm{in}}^2 = 1, \label{innereqn0}
\eeq
with the spike solution 
\beq
y_{\mathrm{in}} = 3 \sqrt{2}\, \sech^2(2^{-1/4}X + A) - 1 - \sqrt{2},
\eeq
where we have matched with the outer solution by requiring that
$y_{\mathrm{in}} \ra - 1 - \sqrt{2}$ as $|X| \ra \infty$.
The constant $A$ (corresponding to a translation in the centre of the
spike) is left undetermined in \cite{bender1999}, although it is
shown in \cite{macgillivray2000} that it must be zero. 
Since this internal spike solution can be combined with any
combination of boundary layers at $x = \pm 1$, we have generated
another four solutions to \eqref{maineqn}.

One might ask whether it is possible to have more than one internal
spike. 
Bender \& Orszag  claim that ``for a given positive value of $\eps$
there are $4(N+1)$ solutions to \eqref{maineqn} which have from 0 to
$N$ internal boundary layers at definite locations, where $N$ is a
finite number depending on $\eps$''.

There seems to have been remarkably little work following up on this
claim. MacGillivray et al.~\cite{macgillivray2000} considered the
solutions with two spikes in detail. They showed that the spikes must be
symmetrically placed about 
$x=0$, and that the separation between them is $O(\eps \log \eps)$. In
view of the rather intricate asymptotic analysis in
\cite{macgillivray2000}, it is perhaps not surprising that no attempt has
been made to analyze the three spike solutions.

On the other hand, the asymptotic dependence of the maximum number of
spikes  $\maxmax$ 
on $\eps$ has been determined. Ai \cite{ai2003} showed that $\maxmax$ is
$O(1/\eps)$, and 
subsequently Wong
and Zhao
\cite{wong2008} showed that
\[ \maxmax \sim \left\lfloor \frac{K}{\eps} \right\rfloor,\]
where $K \approx  0.4725$ and $\lfloor x \rfloor$ is the greatest
integer less than or equal to $x$. Wong
and Zhao also showed that the number of solutions of (\ref{maineqn})
is between $4 \maxmax-3$ and $4 \maxmax$.

We also note that Kath has developed a general method  which
gives a qualitative understanding of the number and type of solutions to
equations such as (\ref{maineqn}) in terms of slowly varying phase
planes \cite{kath1985}. Kath's conclusions are similar to those of
Bender and Orszag. 

In this paper we investigate the claim of Bender \& Orszag, both
numerically and asymptotically. 

We first apply a powerful new algorithm for computing bifurcation
diagrams, deflated continuation \cite{farrell2015d}, to Carrier's
problem. This computation reveals a striking and intricate bifurcation
diagram, with new solutions coming into existence via an apparently infinite
sequence of alternating pitchfork and fold bifurcations as $\eps \ra
0$. Furthermore, its results suggest that the claim of Bender \&
Orszag is incorrect: for each fixed value of $\eps$, the number of
solutions is divisible by 2, but is not always a multiple of 4. (However,
the proportion of values of $\eps$ for which the number of solutions
is not divisible by 4 shrinks rapidly as $\eps \ra 0$.)

We then apply the method of Kuzmak \cite{kuzmak1959} to construct
asymptotic solutions to \eqref{maineqn} with a large number of
internal spikes. This method is a generalisation of both the method of
multiple scales and the WKB method, producing a solution in the form
of a slowly modulated fast oscillation. The frequency and amplitude of
the fast oscillation are allowed to vary slowly with position
(as in the WKB method), but the underlying oscillator is nonlinear, so
that the oscillations are not simply harmonic. We will find that this
asymptotic approach is able to capture very well the bifurcations
identified numerically, providing a more-or-less complete asymptotic
description of the solutions of (\ref{maineqn}).

\section{Numerical analysis and computational results}
\label{num}

\begin{figure}
\centering
\includegraphics[width=12.75cm]{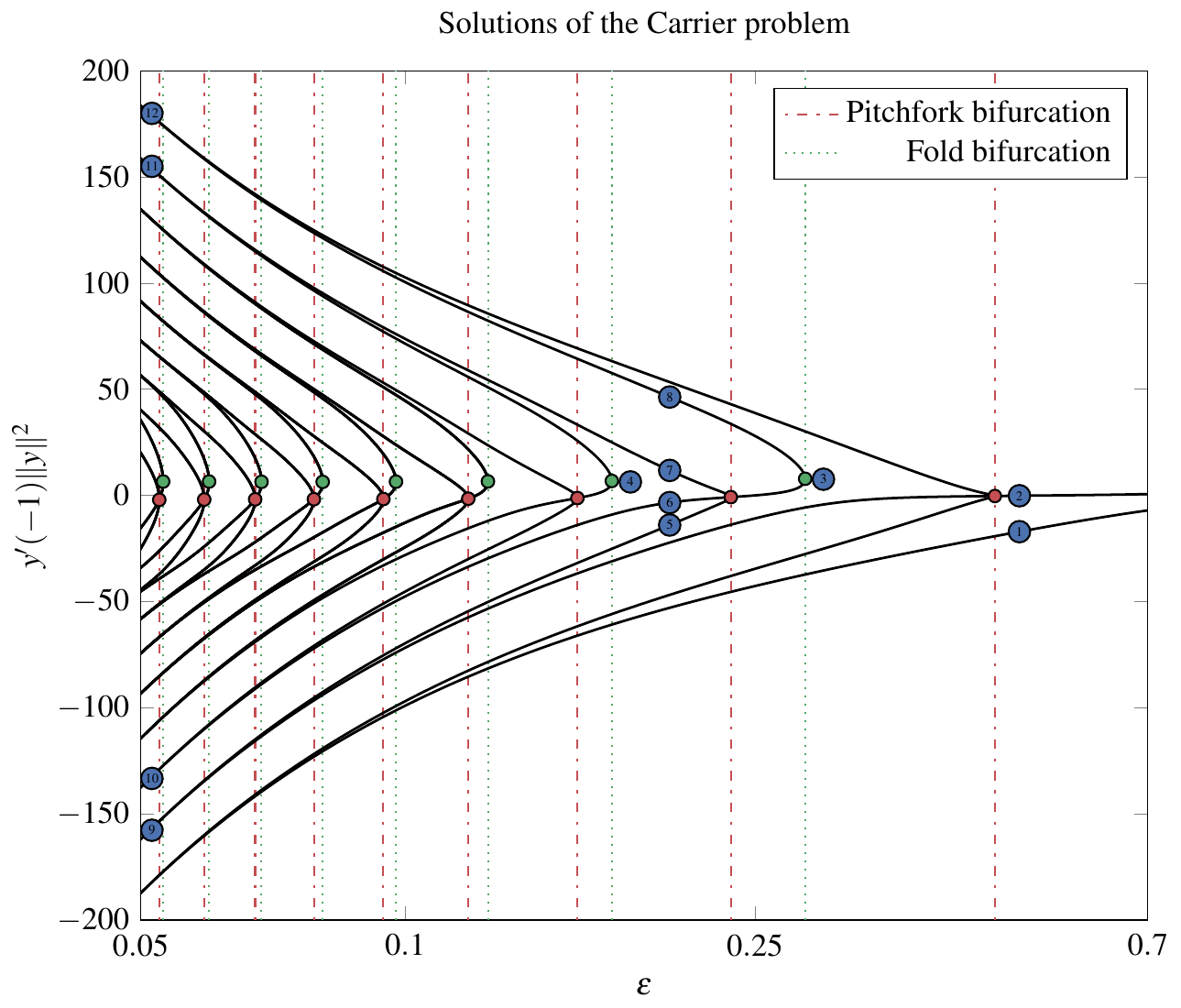}

\begin{tabular}{rrrr}
\includegraphics[width=2.8cm]{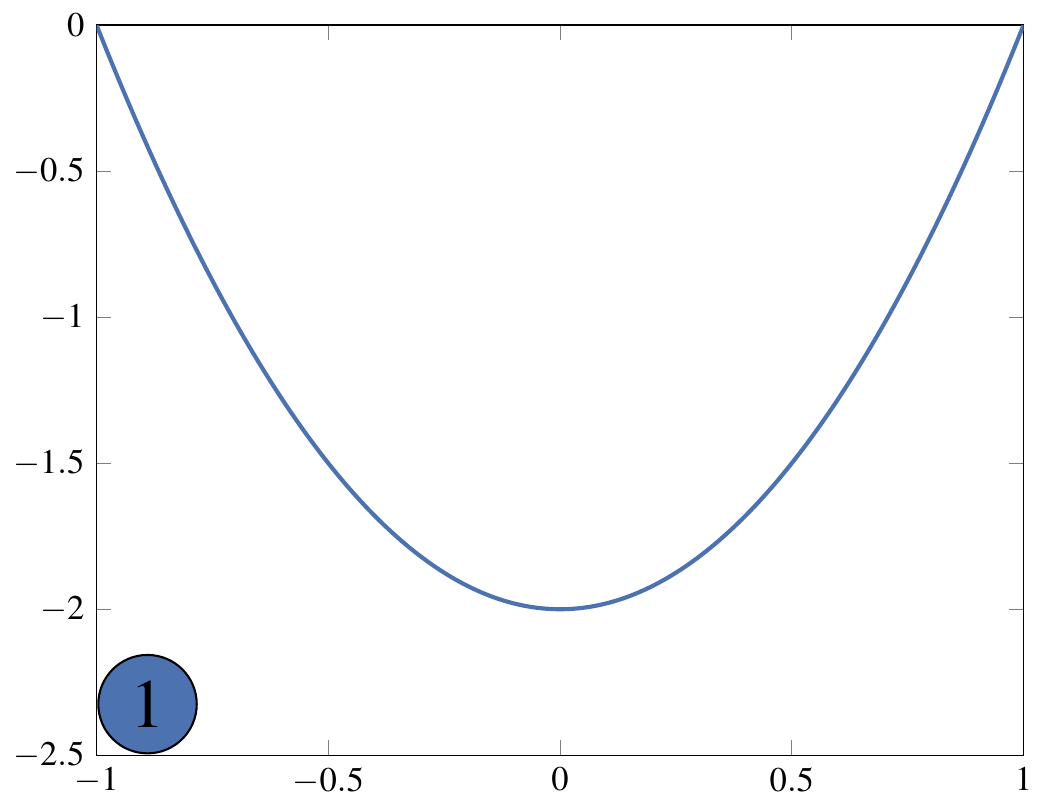} &
\includegraphics[width=2.8cm]{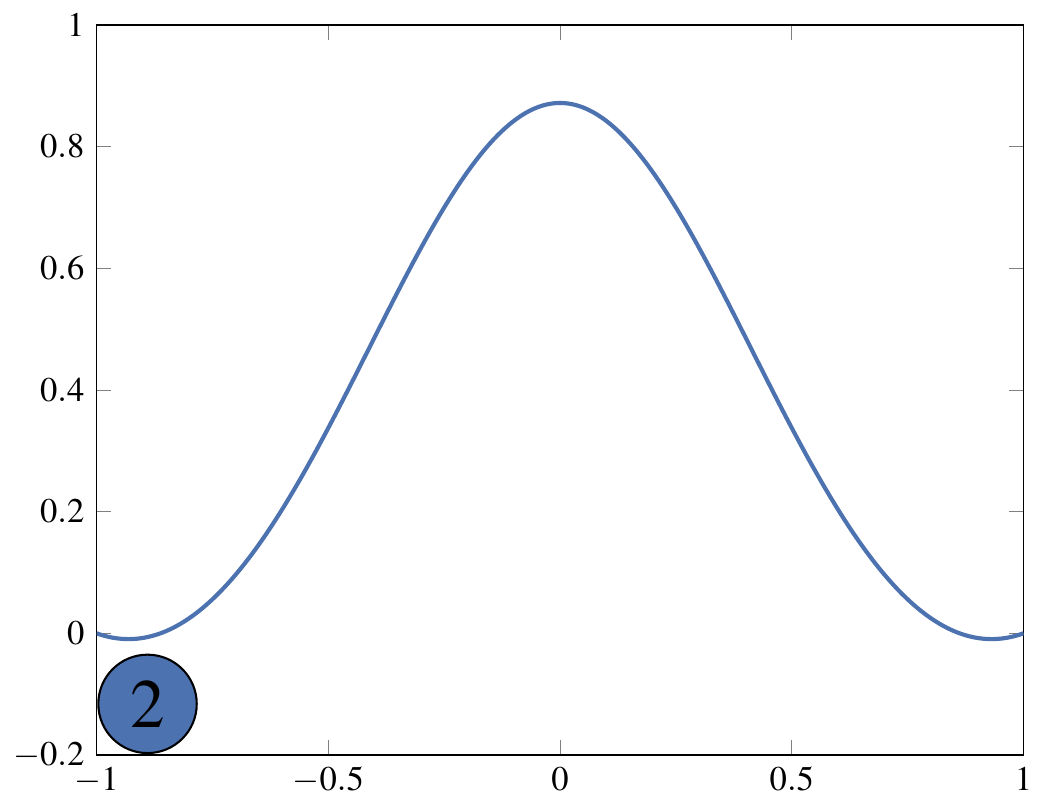} &
\includegraphics[width=2.8cm]{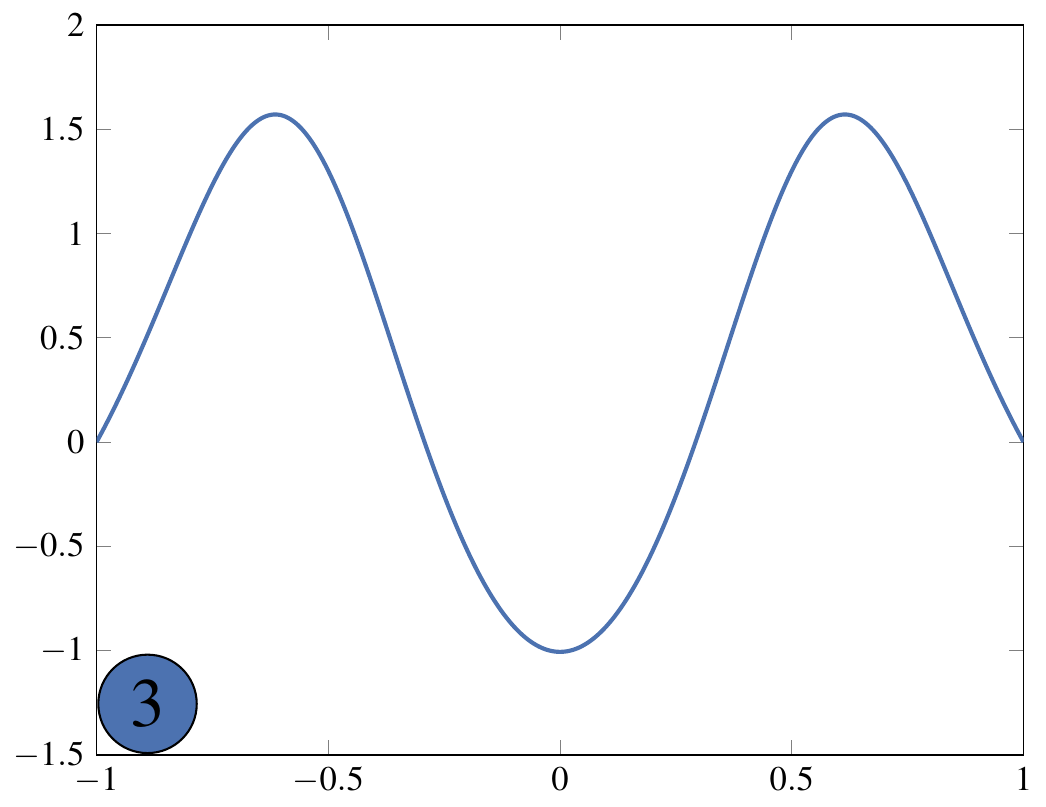} &
\includegraphics[width=2.8cm]{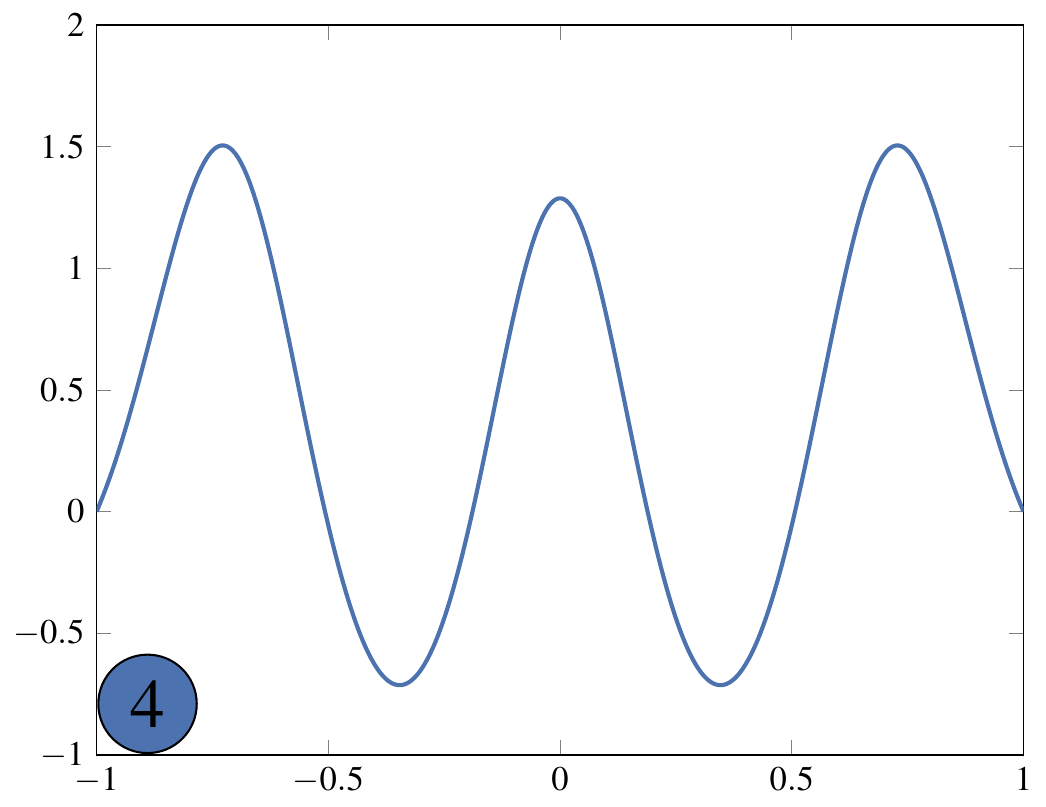}
\\
\includegraphics[width=2.8cm]{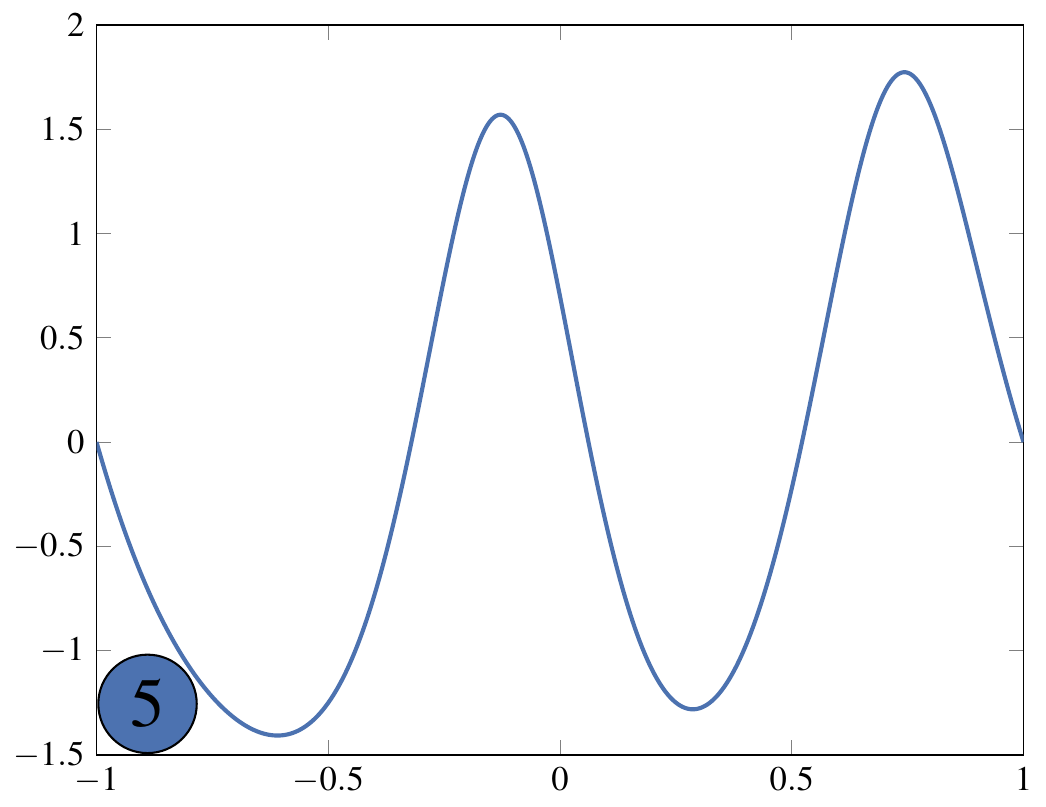} &
\includegraphics[width=2.8cm]{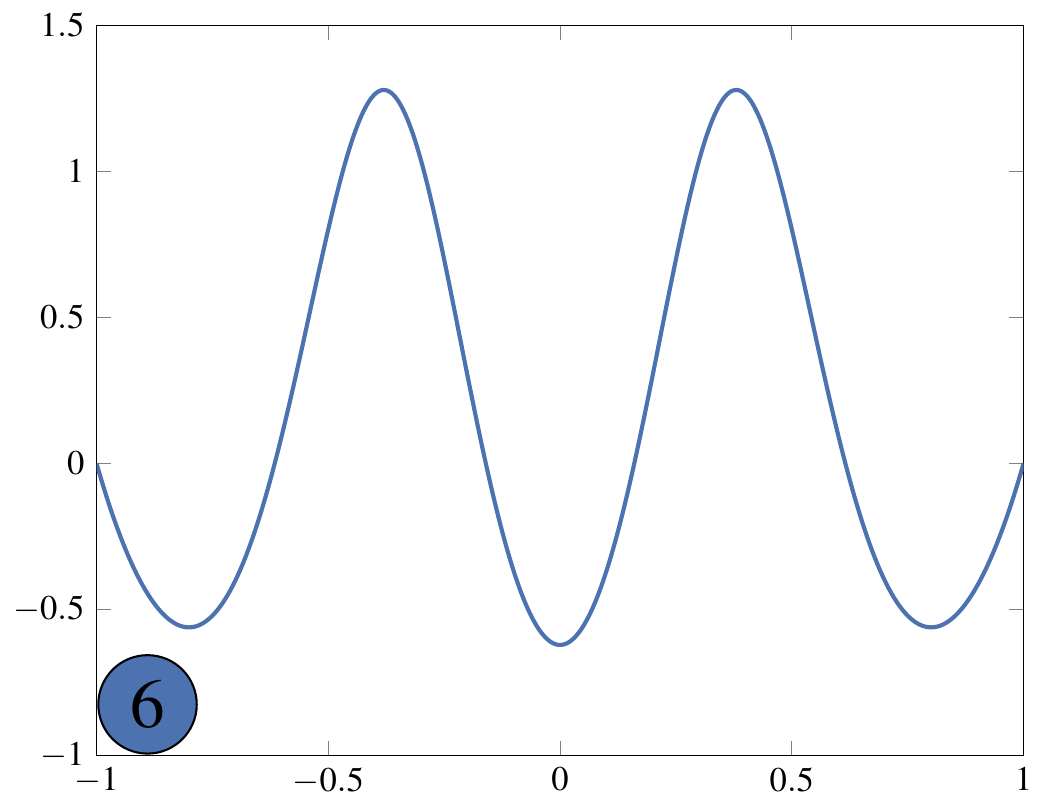} &
\includegraphics[width=2.8cm]{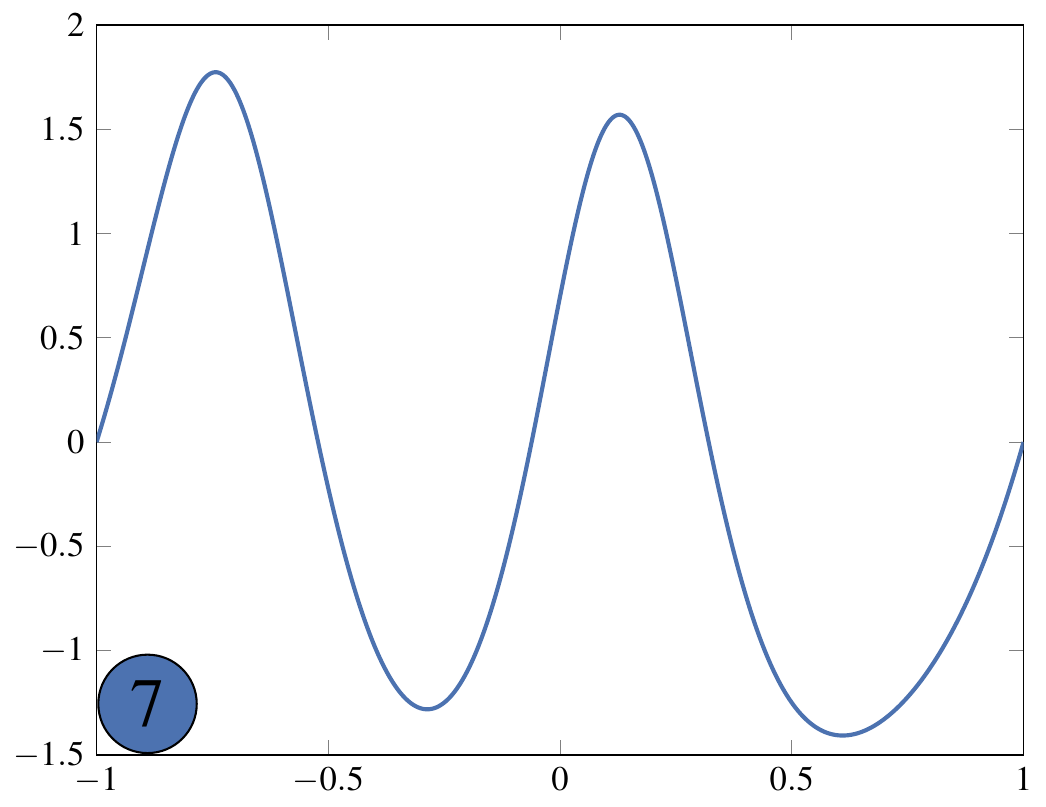} &
\includegraphics[width=2.8cm]{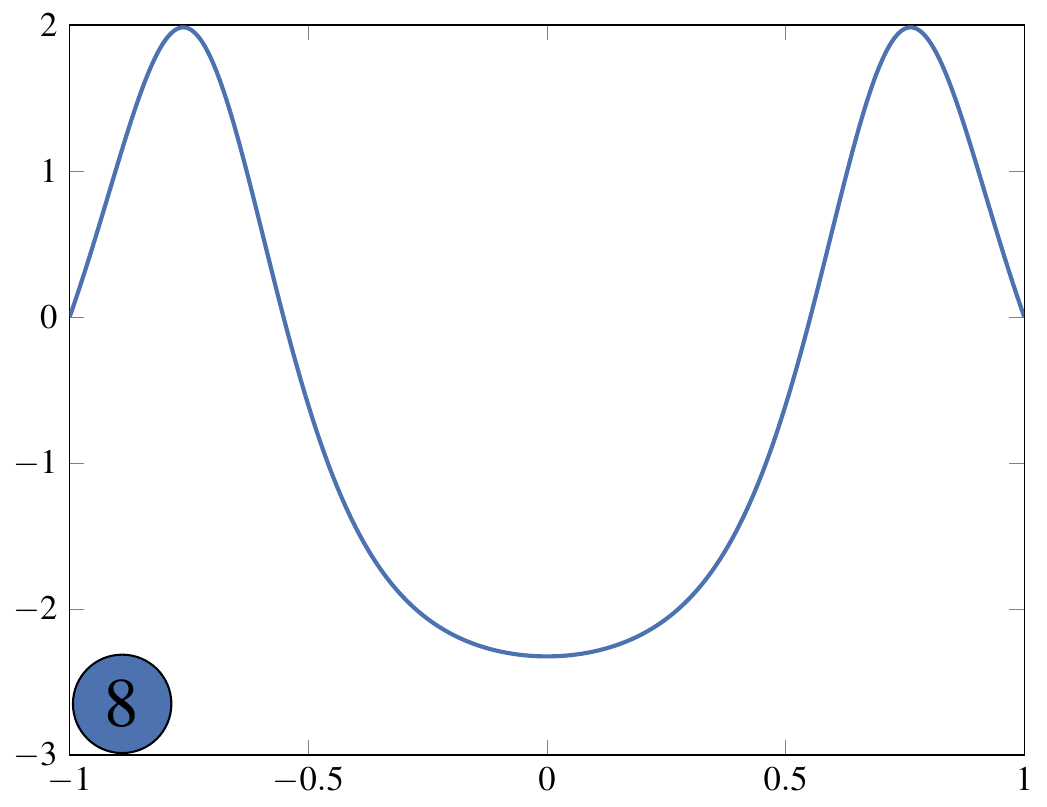}
\\
\includegraphics[width=2.8cm]{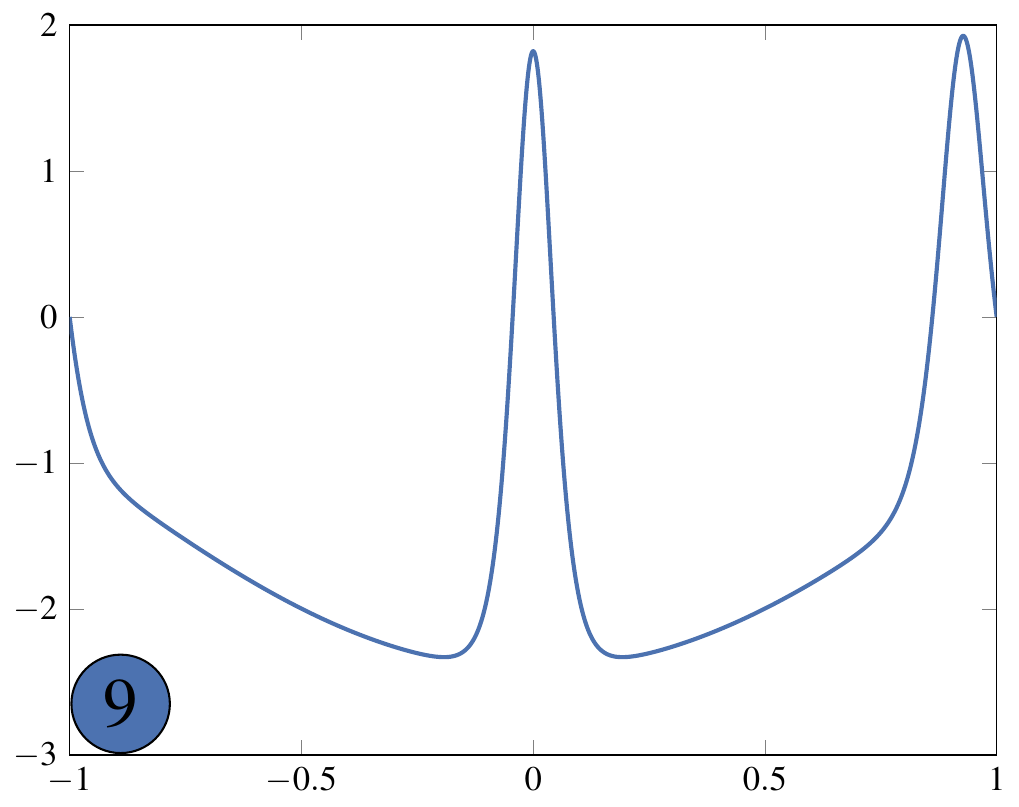} &
\includegraphics[width=2.8cm]{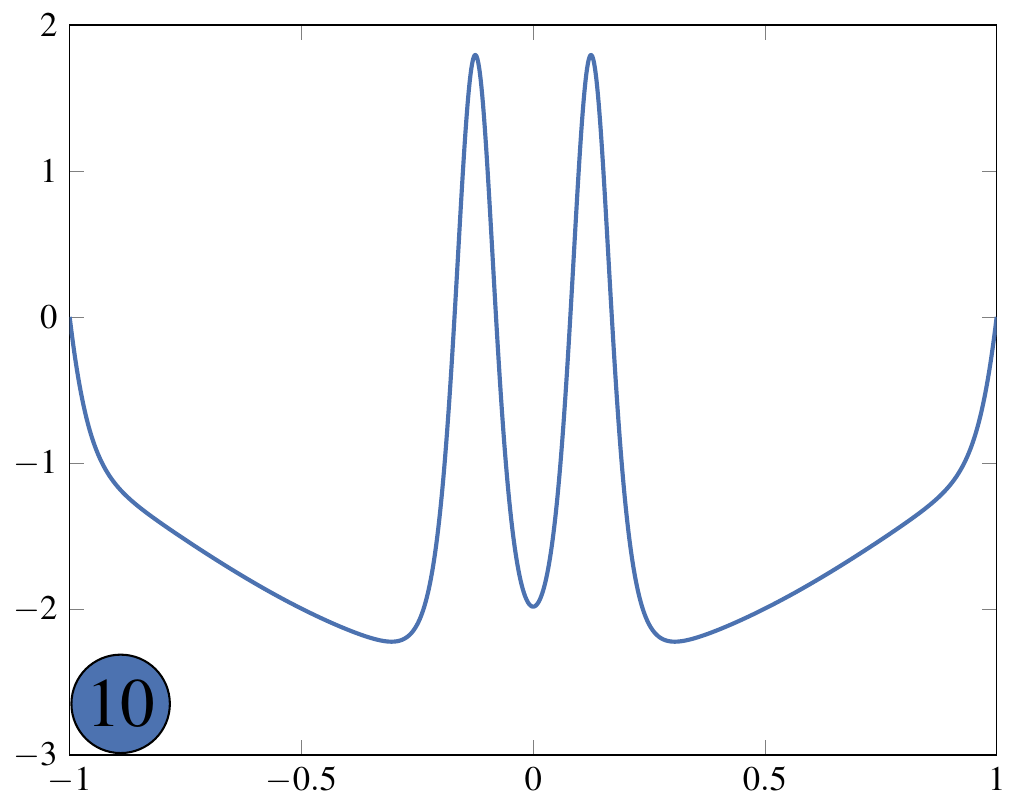} &
\includegraphics[width=2.8cm]{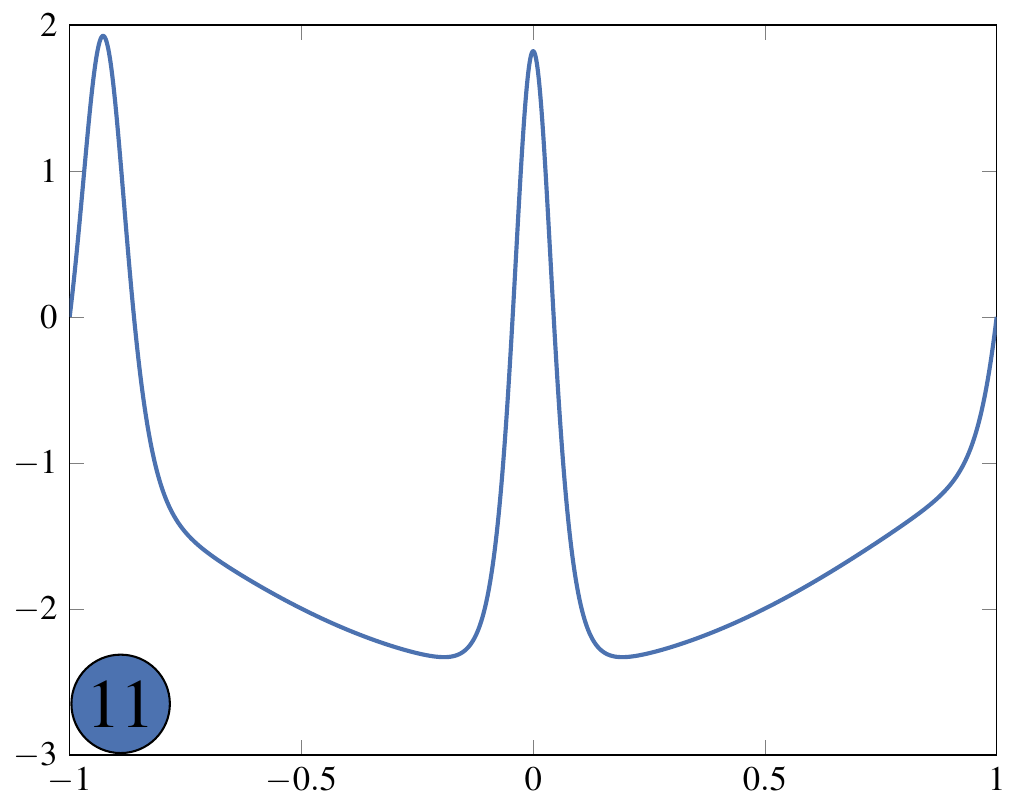} &
\includegraphics[width=2.8cm]{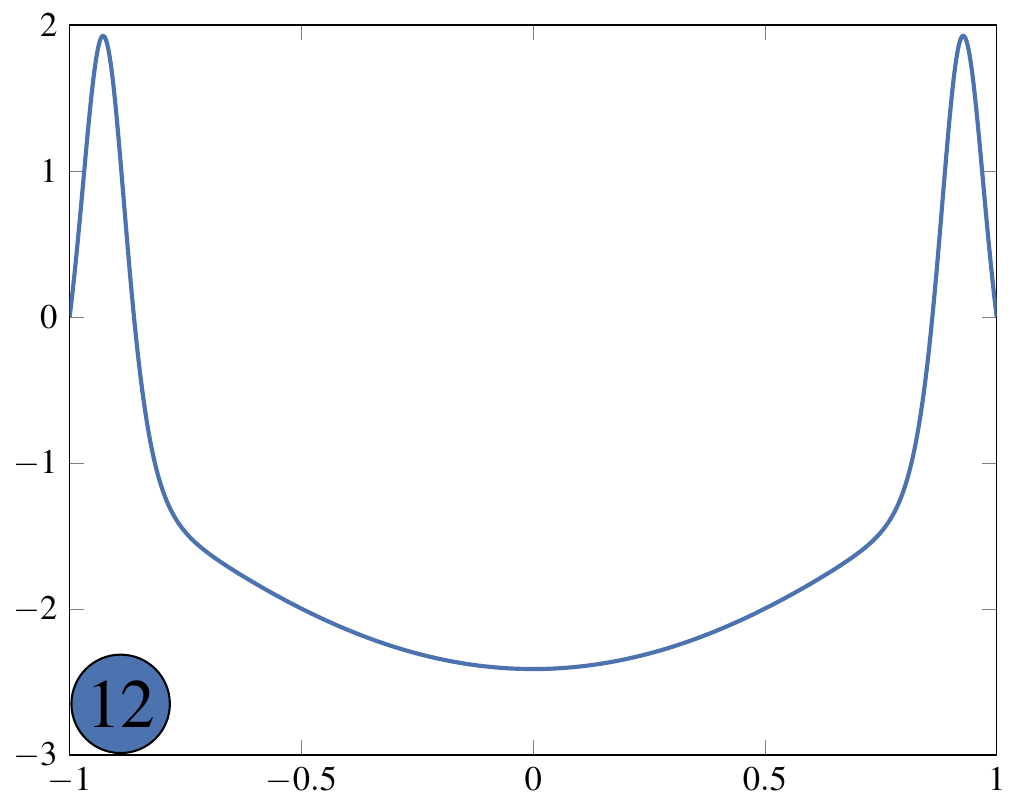}
\end{tabular}

\caption{Bifurcation diagram for Carrier's equation \eqref{maineqn}
as a function of the singular perturbation parameter $\eps$. The diagram
should be read from right to left, as $\eps \ra 0$. For large $\eps$,
there are two solutions; as $\eps \ra 0$, the system undergoes
alternating pitchfork and fold bifurcations (vertical lines). Green circles
denote fold bifurcations; red circles denote pitchfork bifurcations; blue
circles refer to solutions shown in the panels.}
\label{fig:bifurcation-diagram}
\end{figure}

The central task of bifurcation theory is to determine how
the number of solutions to an equation changes as a
parameter is varied.  The main algorithm used to compute
this is the combination of arclength continuation and branch
switching, as invented by Keller in 1977 \cite{keller1977}
and implemented in popular software packages such as AUTO
\cite{doedel1986}. This algorithm is mature and highly
successful, but has a significant drawback: it can only
compute that part of the bifurcation diagram connected to
the initial data, i.e.~it computes connected
\emph{components} of the bifurcation diagram and cannot
``jump'' from one disconnected component to another.
Unfortunately, the bifurcation diagram for Carrier's
problem is indeed disconnected: from any initial datum
branch switching discovers at most four solutions, as we
will see shortly. Since we already know from the analysis
of Section \ref{intro} that at least eight solutions exist
for moderate values of $\eps$, branch switching along $\eps \to 0$ 
offers only a limited insight into the solutions of \eqref{maineqn}.

In recent work, Farrell, Beentjes \& Birkisson have
developed an entirely new algorithm for computing
bifurcation diagrams, called \emph{deflated continuation}
\cite{farrell2015d}. One of the central advantages of
deflated continuation is that it is capable of computing
disconnected bifurcation diagrams, such as that arising in
Carrier's problem\footnote{The other central advantage is
that deflated continuation scales to very large
discretizations, which is more relevant to partial
differential equations than ordinary differential
equations.}. The main weakness of branch switching is that
it relies on identifying critical points at which different
branches meet; this is what renders it incapable of
discovering branches that do not meet the known data at any
point. In contrast, deflated continuation relies instead on
a \emph{deflation} technique for eliminating known solutions
from consideration. Suppose $N$ regular solutions $y_1, y_2,
\dots, y_N$ are known to a discretization of
\eqref{maineqn}. Deflation constructs a new problem residual
for Newton's method with the property that \emph{no initial
guess will converge to $y_1, y_2, \dots, y_N$}.  By
guaranteeing that Newton's method will not converge to known
solutions, deflation enables the discovery of unknown
branches, even if those branches are not connected to the
known data. For more details, see Farrell et al.~\cite{farrell2015d}.

Equation (\ref{maineqn}) was discretized with $5 \times 10^4$ standard
piecewise linear finite elements using FEniCS
\cite{logg2011} and PETSc \cite{balay2015}. We applied
deflated continuation to this discretization from $\eps =
\sqrt{1/2}$ to $\eps = 1/20$, with a continuation step of
$10^{-5}$ in $\eps^2$. Deflation was applied using the $H^1$
norm. All nonlinear systems were solved with Newton's method
and all arising linear systems were solved with LU factorization. Once
deflated continuation had completed, arclength continuation was applied
backwards in $\eps$ from the solutions found at $\eps = 1/20$. The
intricate bifurcation diagram computed in this way is shown
in Figure \ref{fig:bifurcation-diagram}.

The algorithm discovers two solutions to \eqref{maineqn} for $\eps =
\sqrt{1/2}$ from the initial guess ${y}(x) = 1$; at $\eps =
1/20$, 36 solutions were found. The system undergoes an initial pitchfork
bifurcation at $\eps \approx 0.4689$, and subsequently alternates
between fold and pitchfork bifurcations. At each bifurcation, two new
solutions come into existence, and thus there are regions of the
diagram for which the claim of Bender \& Orszag that the number of
solutions is divisible by 4 does not hold. These regions are
precisely the gaps between each fold bifurcation and its
subsequent pitchfork. However, these gaps tend to zero as
$\eps \ra 0$.

The diagram is highly fragmented; no connected component
comprises more than four solutions, which is why branch
switching applied to this problem can never discover more
than four solutions from any single initial datum. We observe
that each connected component is characterized by the
number of interior maxima $\nummax$. 
The lowest component comprises
a single branch that undergoes no bifurcations, and
corresponds to $\nummax = 0$; this branch exists for all $\eps$ and
is shown in panel 1 of Figure \ref{fig:bifurcation-diagram}.
The next component ($\nummax = 1$) also exists for all $\eps$, and
is shown in panel 2 of Figure \ref{fig:bifurcation-diagram}.
After the  pitchfork bifurcation it comprises three
solutions, one with an interior spike but no boundary spikes, 
one with no interior spike but a boundary spike on
the left (as in Figure \ref{fig1}(b)), and one with no interior spike but
a boundary spike on 
the right (as in Figure \ref{fig1}(c)).

\begin{figure}
\centering
\subfigure[]{
\scalebox{0.5}{\includegraphics{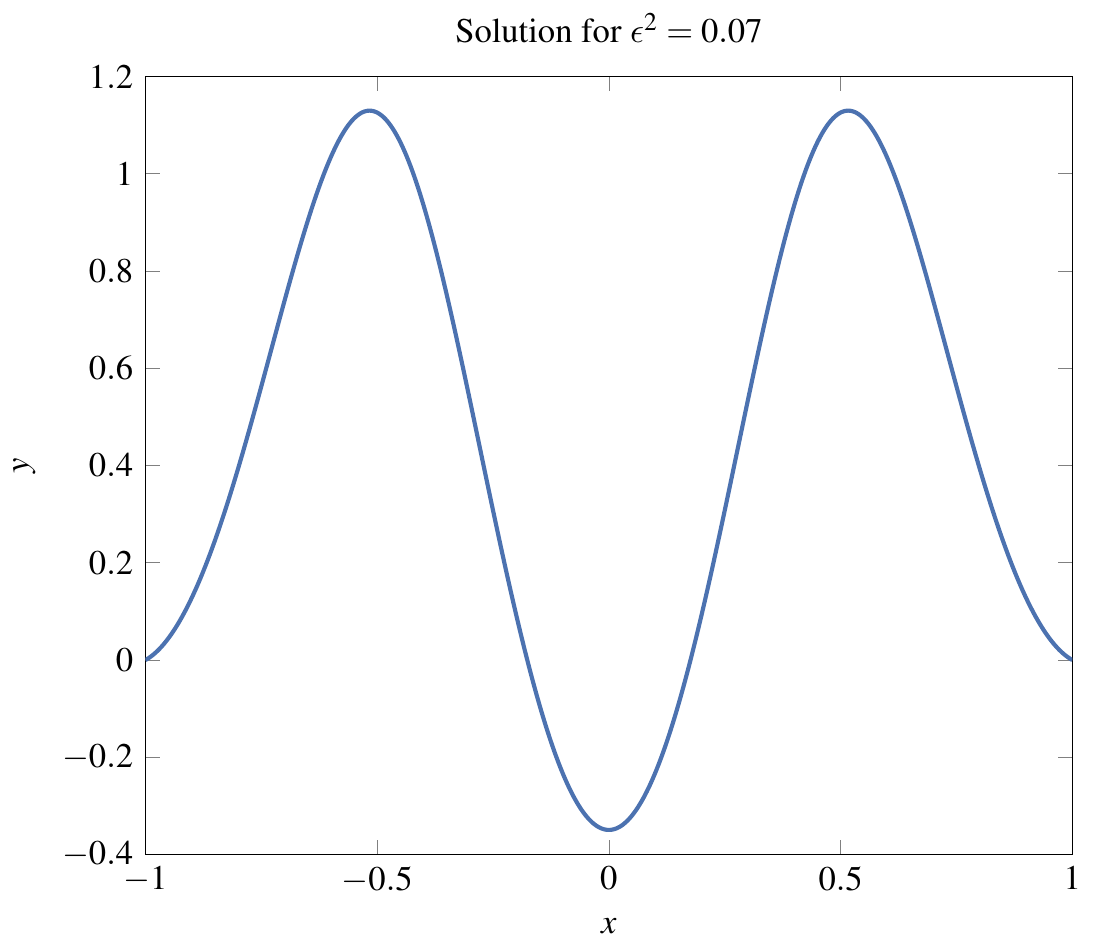}}
}
\
\subfigure[]{
\scalebox{0.5}{\includegraphics{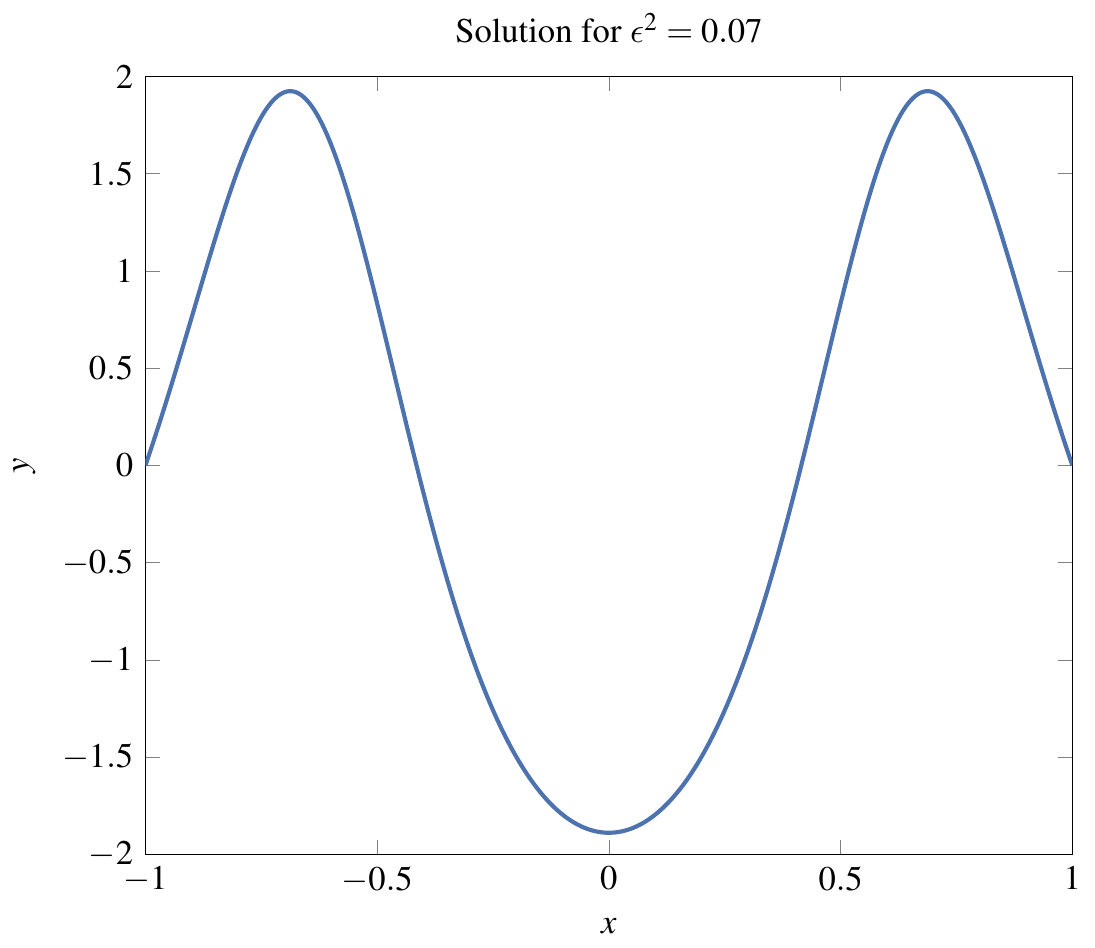}}
}
\
\caption{The two solutions for $\nummax = 2$ at $\eps^2 =
0.07$, between the first fold (at $\eps^2 \approx 0.08135$)
and its subsequent pitchfork (at $\eps^2 \approx 0.05509$).
The maxima in (a) are more concentrated towards the centre,
whereas (b) has maxima near the boundaries.}
\label{fig:epssq-0.07}
\end{figure}

\begin{figure}
\centering
\subfigure[]{
\scalebox{0.5}{\includegraphics{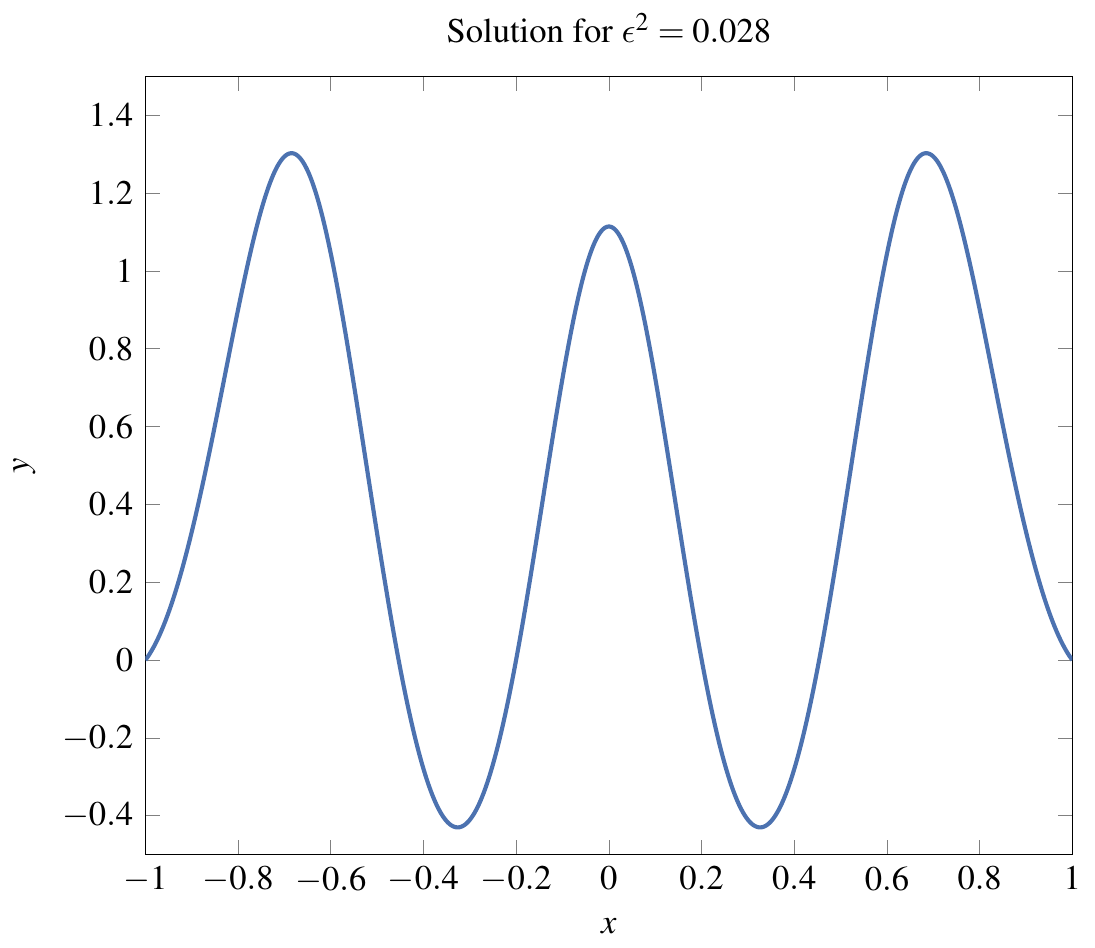}}
}
\
\subfigure[]{
\scalebox{0.5}{\includegraphics{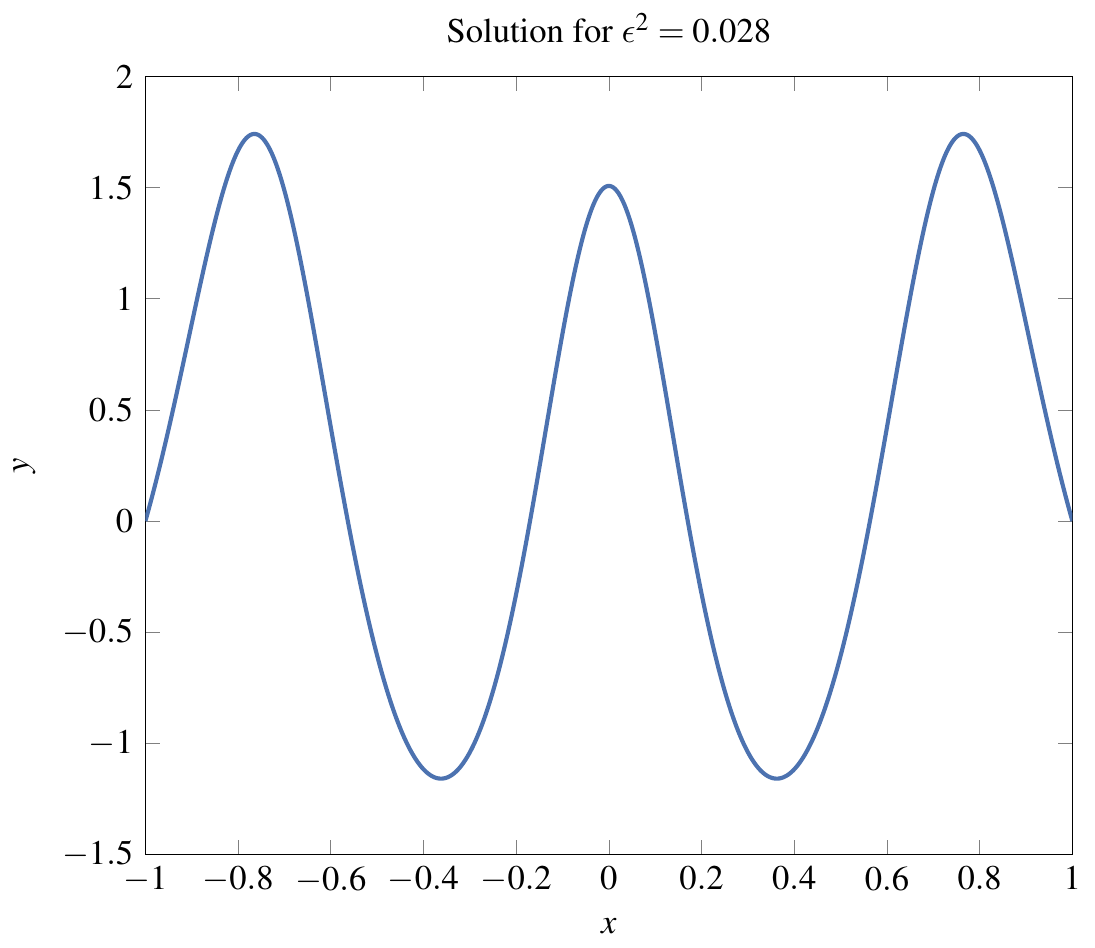}}
}
\
\caption{The two solutions for $\nummax = 3$ at $\eps^2 =
0.028$, between its originating fold (at $\eps^2 \approx 0.02953$)
and its subsequent pitchfork (at $\eps^2 \approx 0.02466$).
The maxima in (a) are more concentrated towards the centre,
whereas (b) has maxima near the boundaries. Compare to Figure
\ref{fig:epssq-0.07}, the corresponding diagram for $\nummax = 2$.}
\label{fig:epssq-0.028}

\subfigure[]{
\scalebox{0.5}{\includegraphics{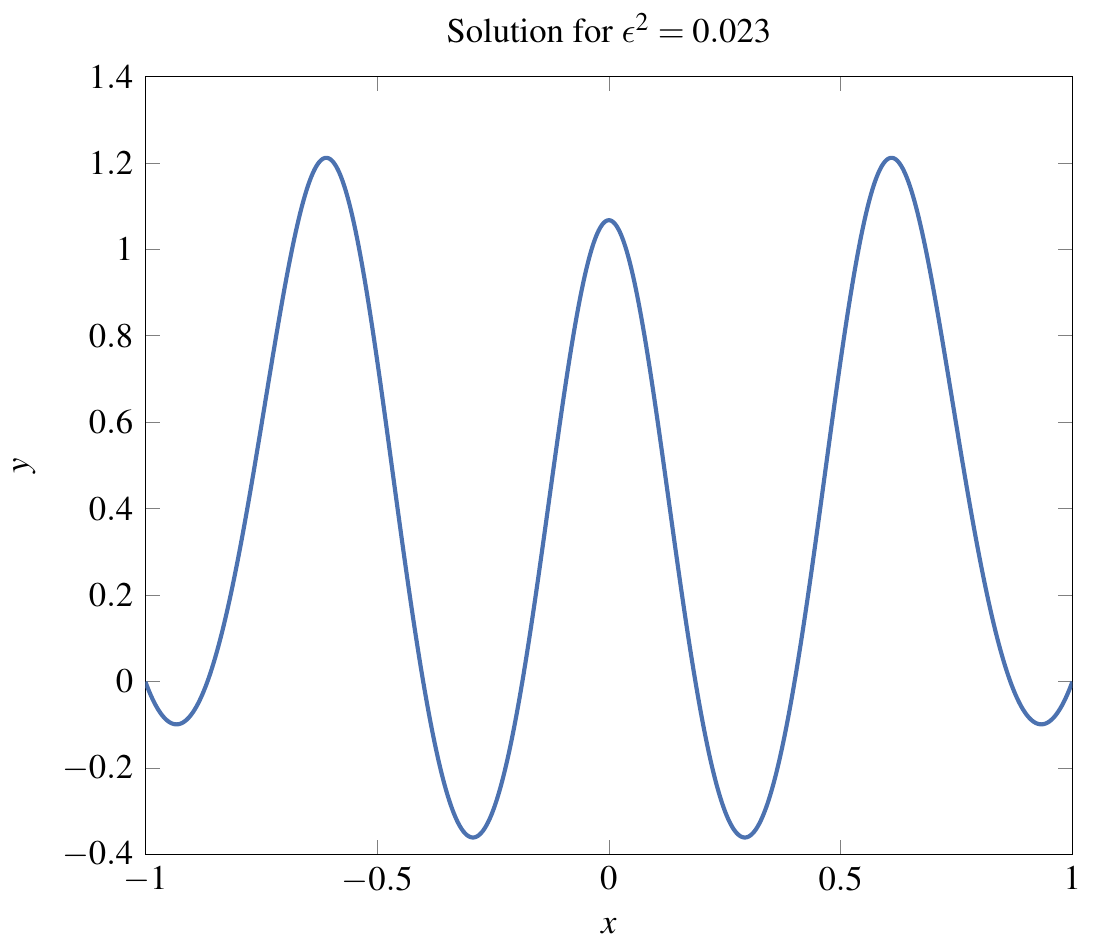}}
}
\
\subfigure[]{
\scalebox{0.5}{\includegraphics{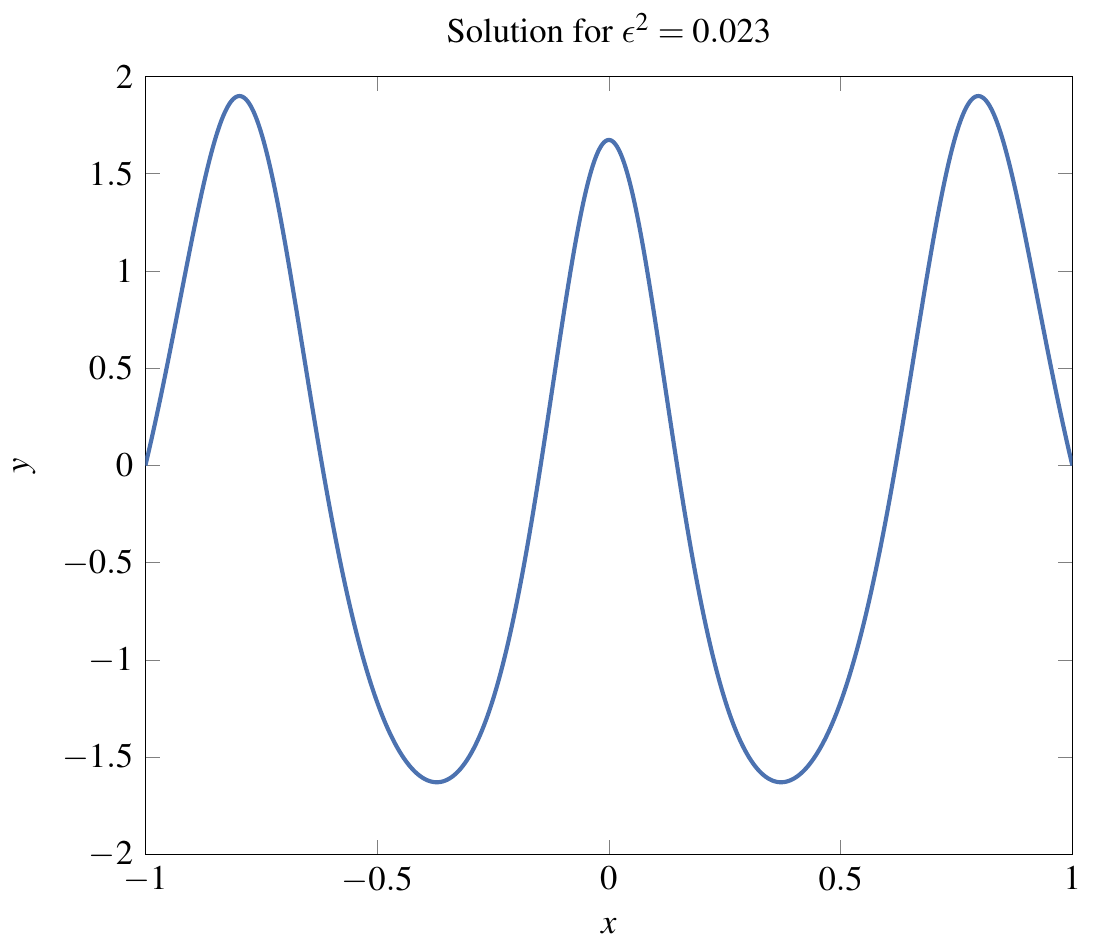}}
}
\\
\subfigure[]{
\scalebox{0.5}{\includegraphics{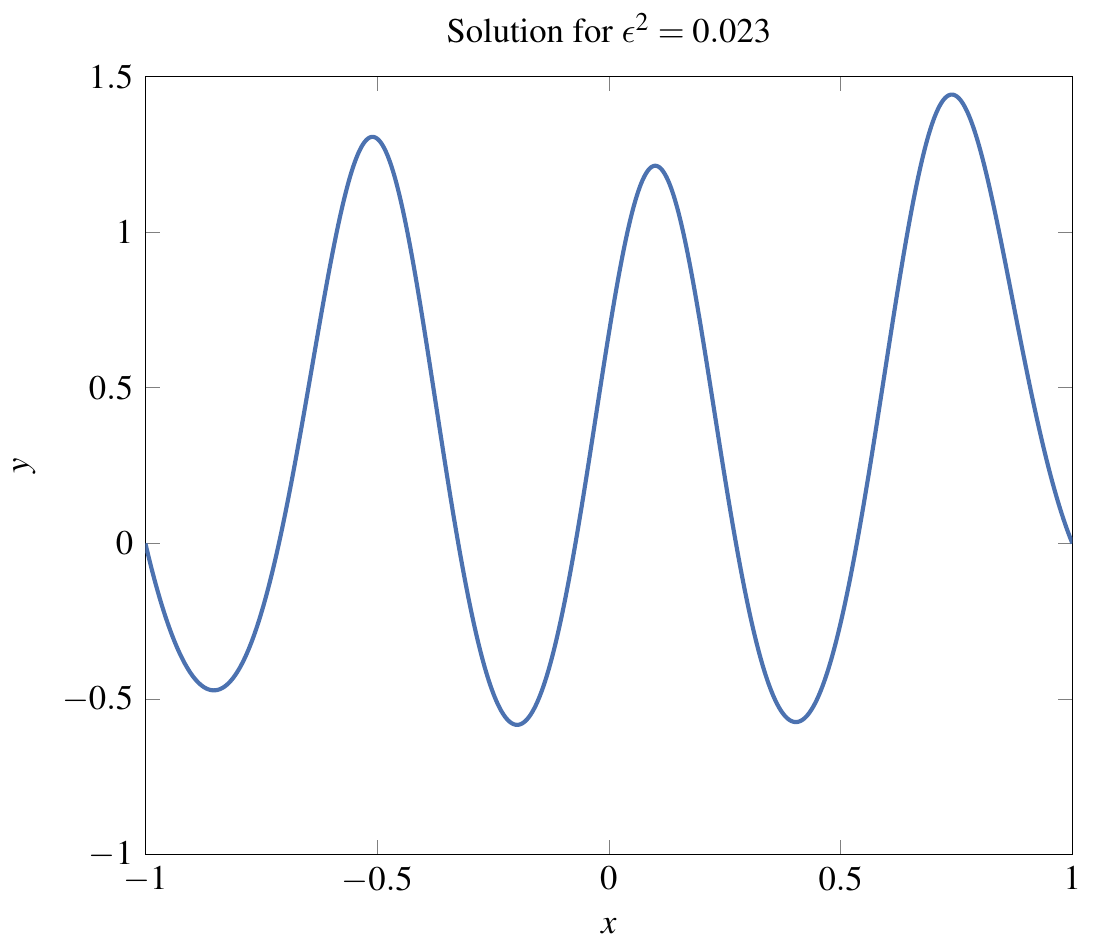}}
}
\
\subfigure[]{
\scalebox{0.5}{\includegraphics{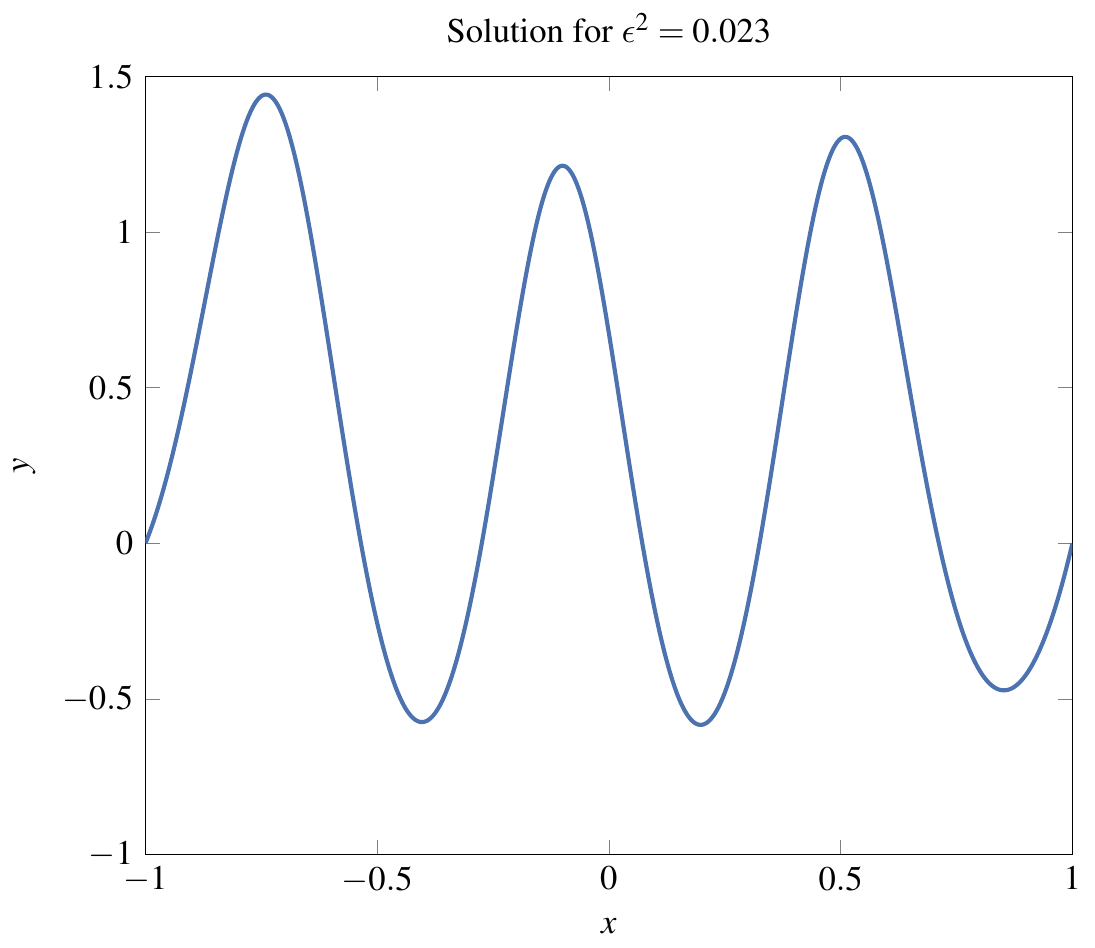}}
}
\caption{The four solutions for $\nummax = 3$ at $\eps^2 = 0.023$,
after its pitchfork bifurcation at $\eps^2 \approx 0.02466$.
The two symmetric solutions (a) and (b) lie on the same
branch as \ref{fig:epssq-0.028}a and \ref{fig:epssq-0.028}b
respectively. The branches (c) and (d) have bifurcated from
(a) and are characterized by one of the interior maxima
approaching the boundary. Compare to the middle row of panels of Figure
\ref{fig:bifurcation-diagram}, the corresponding solutions for $\nummax =
2$.}
\label{fig:epssq-0.023}
\end{figure}

\begin{figure}
\centering
\subfigure[]{
\scalebox{0.5}{\includegraphics{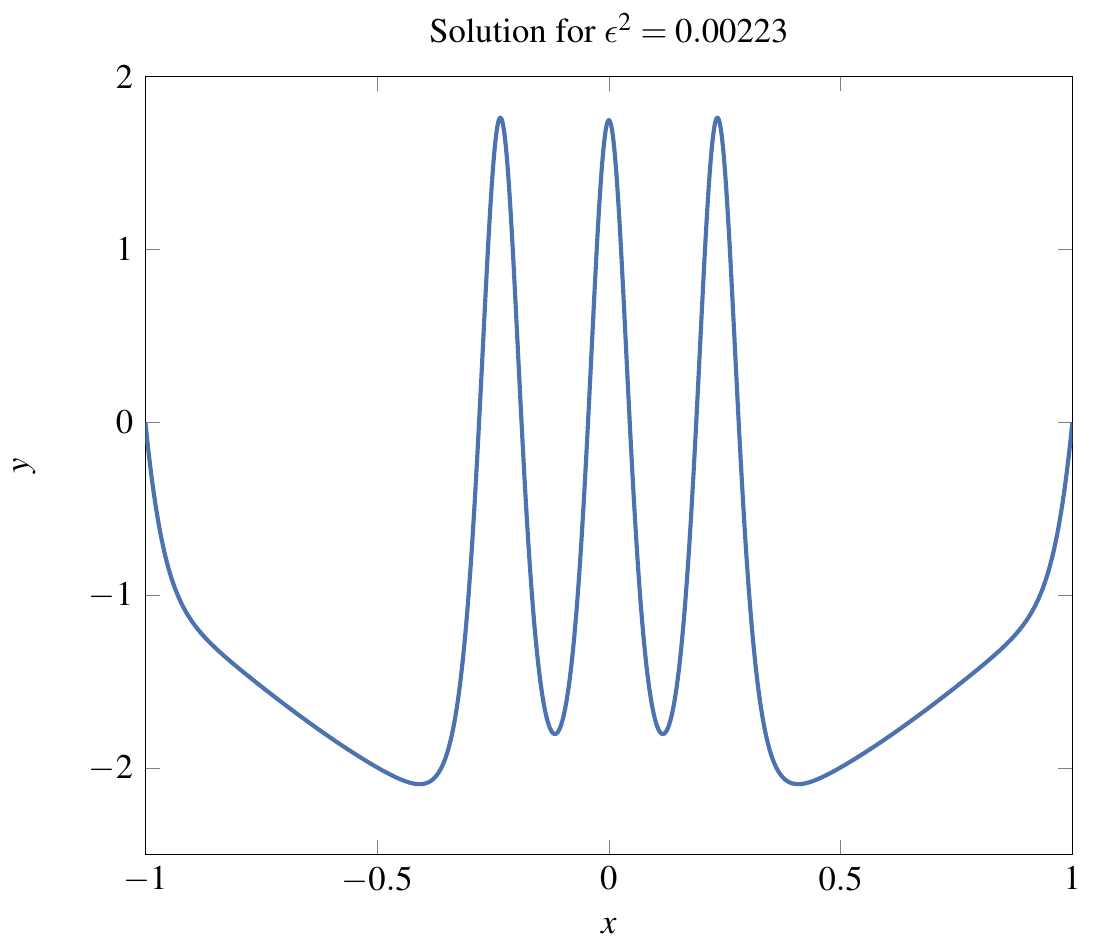}}
}
\
\subfigure[]{
\scalebox{0.5}{\includegraphics{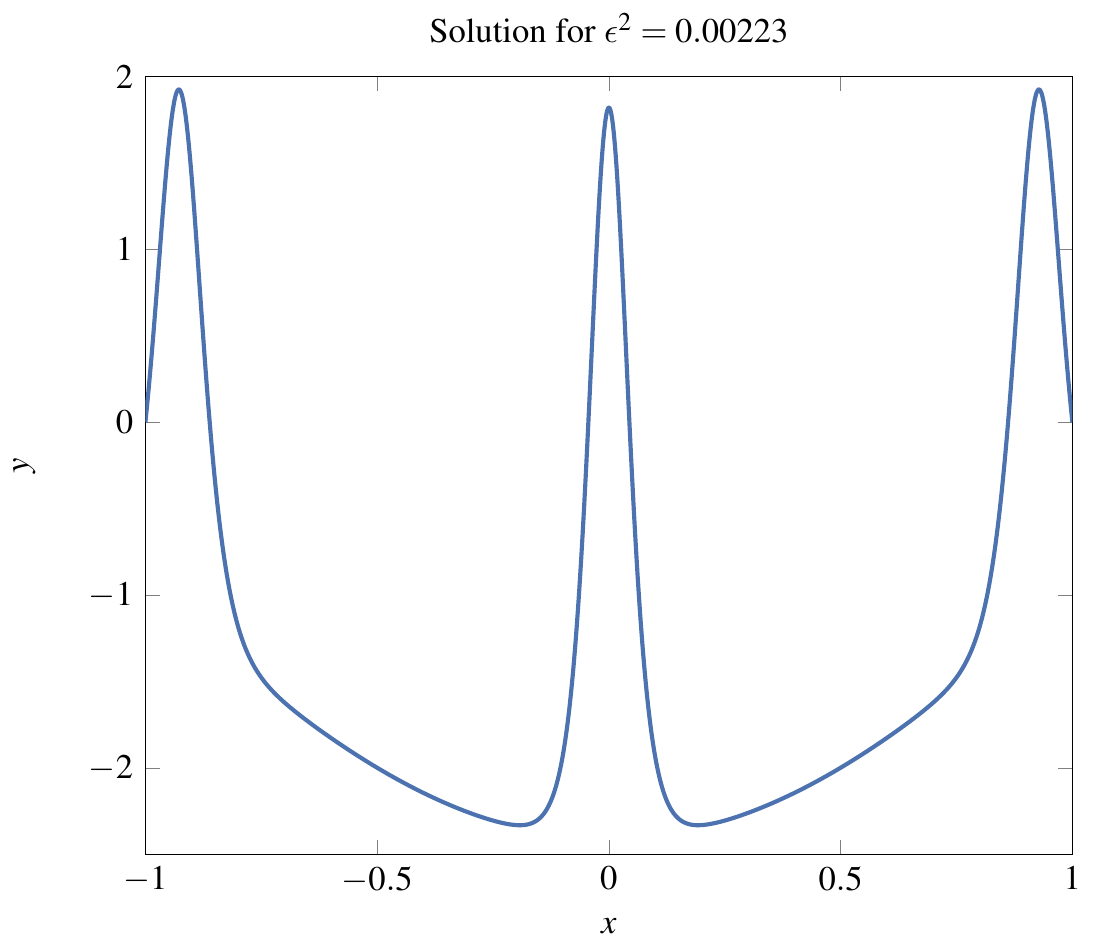}}
}
\\
\subfigure[]{
\scalebox{0.5}{\includegraphics{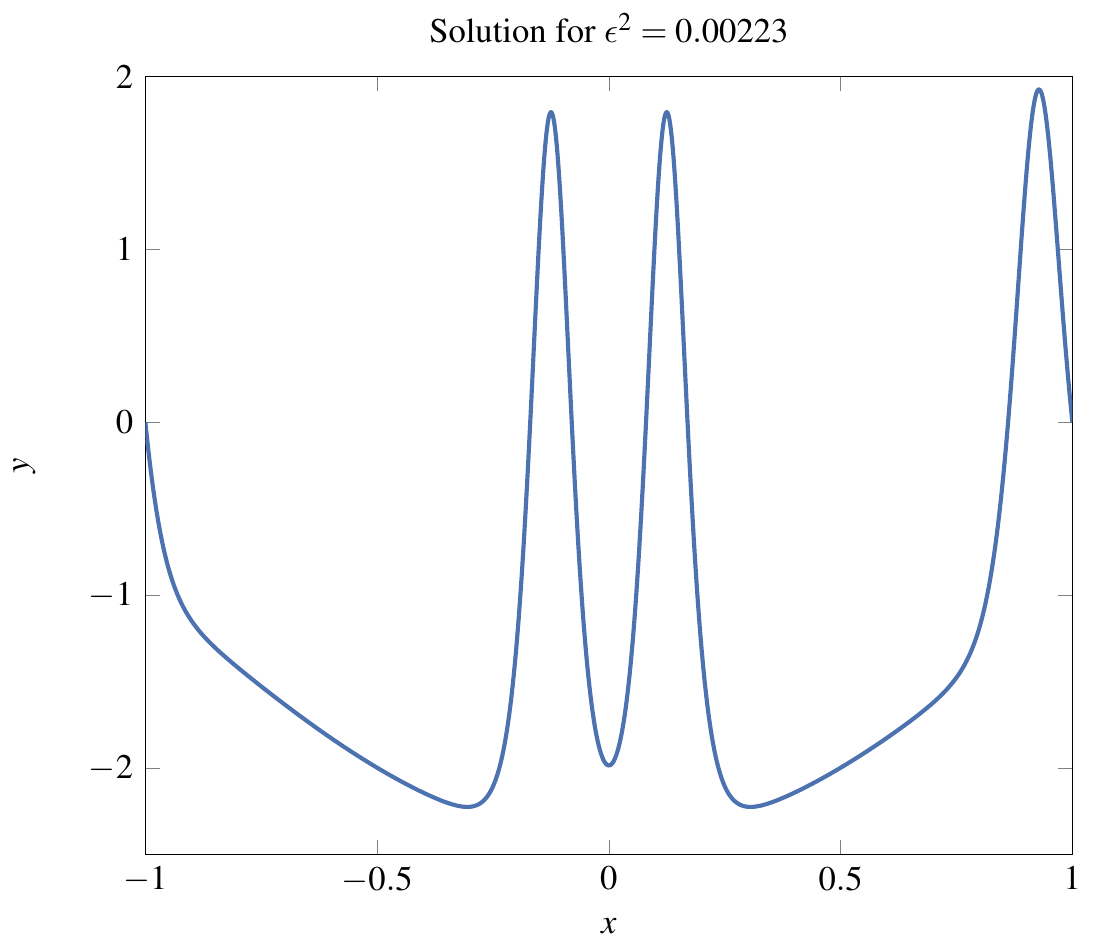}}
}
\
\subfigure[]{
\scalebox{0.5}{\includegraphics{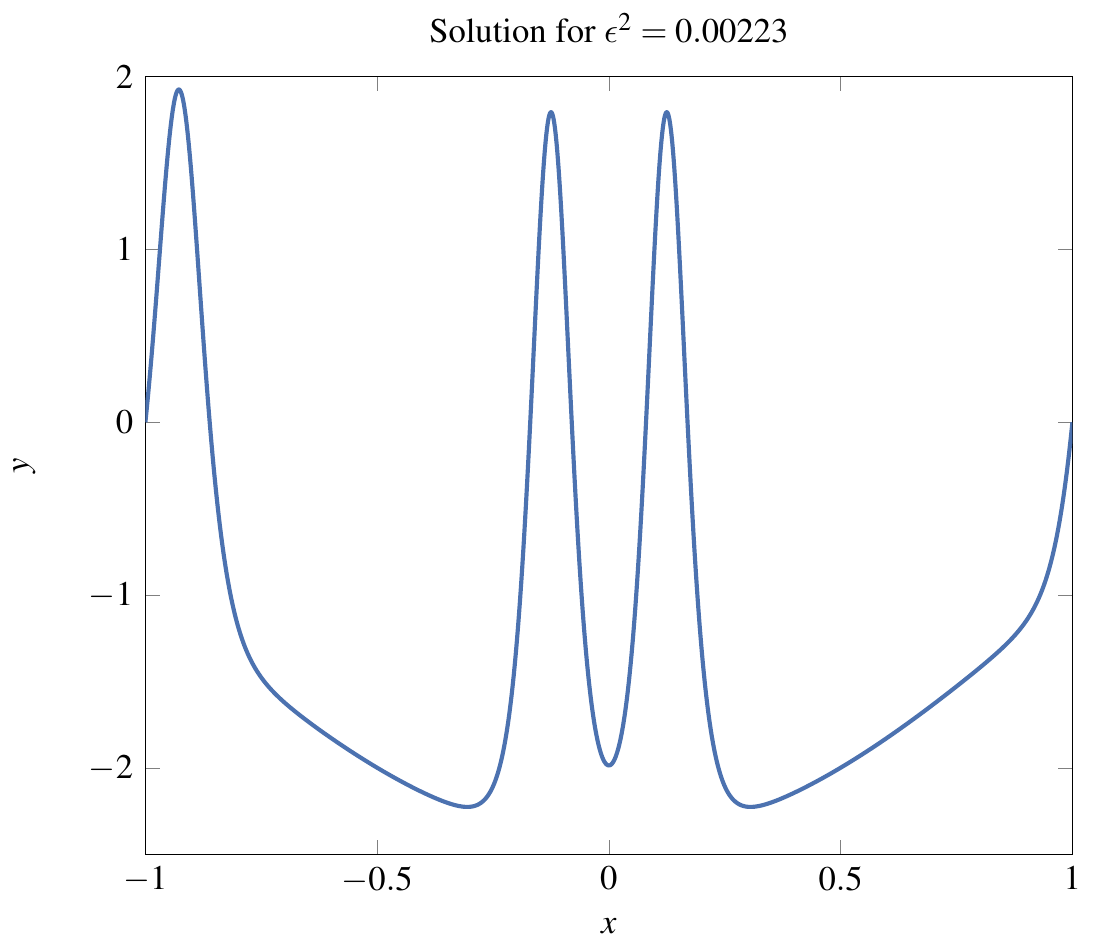}}
}
\caption{The four solutions for $\nummax = 3$ at $\eps^2 =
0.00223$. The solutions are labeled as in
\ref{fig:epssq-0.023}. The symmetric solution with maxima
near the centre has $\nummax$ interior layers; the symmetric
solution with maxima near the boundary has 2 boundary layers
and $\nummax - 2$ interior layers; the asymmetric branches have 1
boundary layer and $\nummax - 1$ interior layers. Compare to
the bottom row of panels of Figure \ref{fig:bifurcation-diagram}, the corresponding diagram
for $\nummax = 2$.}
\label{fig:epssq-0.00223-2}
\end{figure}

The other components do not exist for large $\eps$, and
come into existence at fold bifurcations as $\eps$ is reduced (panels
3 and 4 of Figure \ref{fig:bifurcation-diagram}).
Between each fold and its subsequent pitchfork bifurcation,
the two solutions are symmetric with $n$ local maxima; both
begin with $y'(-1) > 0$. In one of these solutions the
maxima are more concentrated near $x = 0$, whereas in the
other there are maxima near $x = \pm 1$. This is illustrated
in Figure \ref{fig:epssq-0.07} for $\nummax = 2$ and
Figure \ref{fig:epssq-0.028} for $\nummax = 3$. The
symmetry-breaking pitchfork bifurcation occurs on the branch
where the maxima are concentrated near $x = 0$ close to (but
not exactly at) the value of $\eps$ for which $y'(\pm 1) =
0$. The new branches consist of solutions where one maximum
leaves the centre and approaches one of the boundaries. This is
illustrated in the middle row of panels of Figure \ref{fig:bifurcation-diagram} for $\nummax =
2$ and Figure \ref{fig:epssq-0.023} for $\nummax = 3$. As
$\eps \ra 0$, the symmetric solution with maxima concentrated
near the centre tends to a solution with $\nummax$ interior spikes; the
symmetric solution with maxima near $x = \pm 1$ tends to a
solution with 2 boundary layer spikes and $\nummax - 2$ interior spikes;
the solutions emanating from the pitchfork bifurcation tend
to solutions with 1 boundary layer spike and $\nummax - 1$ interior
spikes. This is illustrated in the bottom row of panels of Figure
\ref{fig:bifurcation-diagram} for $\nummax = 2$ and Figure
\ref{fig:epssq-0.00223-2} for $\nummax = 3$.

There are two further remarks to make regarding these
observations. The first is that the number of interior
maxima $\nummax$  does not change through the bifurcations;
hence our observation that each disconnected component is
characterized by $\nummax$. One consequence of this is that
the four solutions we can generate  by choosing different
combinations of boundary layers (\ref{inner}) for a given
number of interior spikes are not all on the same component
of the bifurcation diagram. For example, in Figure
\ref{fig1}, solutions (b) and (c) are connected in the
bifurcation diagram, but lie on a different component to
solutions (a) and (d), which are themselves on different
components.

The second is that close to the bifurcations the
oscillations fill the domain. Thus a boundary layer analysis
in which there are a finite number of interior spikes
separated from boundary layers by a spike-free outer region
will never be able to capture the bifurcations. To capture
the bifurcations we need to generate asymptotic solutions in
which the interior spikes go all the way to the boundary.

Deflated continuation has successfully revealed an enormous
amount of information regarding the solutions to
\eqref{maineqn}, but it does not explain why the bifurcation
diagram possesses this structure or predict the locations of
the alternating fold and pitchfork bifurcations.  Using the
intuition we have gained from the numerical results, we now
turn to asymptotic methods to see if we can predict these
features analytically.

\section{Asymptotic analysis}
\subsection{Asymptotic approximation using Kuzmak's method}
\label{kuzmak}
To be able to predict the bifurcations we have seen in solutions of
\eqref{maineqn} we need to be able to generate asymptotic solutions in
which the interior spikes fill the domain, that is, solutions which are
rapidly oscillating.
 We  construct such 
asymptotic solutions using the 
method of Kuzmak \cite{kuzmak1959}.

We need to allow the frequency of the oscillation to vary slowly. Thus
(as in the WKB method) we define the fast scale as $X = \phi(x)/\eps$,
where the function 
$\phi(x)$ is to be determined. We then look for solutions
$y(x,X)$, treating the slow scale $x$ and the fast scale $X$ as
independent. We remove the indeterminacy this generates, and avoid
secular terms in $X$,  by imposing
exact periodicity in $X$ with period 1.
 
From the chain rule we have
\beqas
\fdd{y}{x} & = & y_x + \frac{\phi'}{\eps} y_X,\\
\sdd{y}{x} & = & y_{xx} +  \frac{2\phi'}{\eps} y_{xX} + \frac{\phi''}{\eps}
y_X + \frac{(\phi')^2}{\eps^2} y_{XX},
\eeqas
where $\phi' = \fdd{\phi}{x}$, and a subscript denotes partial differentiation.
Thus  equation \eqref{maineqn} becomes
\beq
 (\phi')^2 y_{XX} + \eps(2 \phi' y_{xX} + \phi'' y_{X}) + \eps^2
y_{xx} + 2(1-x^2) y + y^2 = 1.\label{eqn}
\eeq
We now pose a series expansion  in powers of $\eps$:
\[ y \sim y_0 + \eps y_1 + \cdots.\]
At leading order this gives
\beq
(\phi')^2 \spd{y_0}{X} + 2(1-x^2) y_0 + y_0^2 = 1,\label{leading}
\eeq
with $y_0$ periodic in $X$, with period 1.
If (\ref{leading}) were a linear equation then the solutions would be
exponentials and our asymptotic method would simply be the WKB
method. The fact that (\ref{leading}) is nonlinear means the fast
oscillator is not simply harmonic, and we have to work a little bit
harder to describe the oscillations.

Multiplying (\ref{leading}) by $2\pd{y_0}{X}$ and integrating gives
\[(\phi')^2 \left(\pd{y_0}{X}\right)^2 + 2(1-x^2) y_0^2 + \frac{2}{3}y_0^3 =
2y_0 + A(x),\]  
where the constant of integration, $A$, depends on the slow scale $x$.
Separating the variables and integrating again gives
\beq
 \pm\phi' \int_0^{y_0}\frac{\d y}{(A(x) + 2y - 2(1-x^2)y^2 - 2
  y^3/3)^{1/2}}  = X + \mu(x) .\label{1}
\eeq
The relevant values of $A(x)$ are those for which the cubic in the
denominator has three 
real roots, $Y_0 <Y_1 <Y_2$, say. The cubic is positive for values
$Y_1<y<Y_2$. Periodicity of $y_0$ in $X$ is achieved by integrating
the positive square root from $Y_1$ to $Y_2$ and then the negative
square root from 
$Y_2$ to $Y_1$. Setting the period equal to unity therefore implies
\beq
\phi' = \left(2\int_{Y_1}^{Y_2} \frac{\d y}{(A(x) + 2y - 2(1-x^2)y^2 - 2
  y^3/3)^{1/2}} \right)^{-1}.\label{phiprime}
\eeq
Note\footnote{In the WKB method
$\phi'$ would be independent of $A$.} that $\phi'$ depends on $A$ but not on $\mu$.

We have thus far determined the leading-order solution $y_0$ and the fast
scale of oscillation $\phi$ in terms of the two unknown slow functions $A(x)$
and $\mu(x)$, which correspond roughly to the slowly modulated
amplitude and  phase of the oscillation. To determine these functions
we need to proceed to higher orders in the asymptotic expansion.

Equating coefficients of $\eps$ in (\ref{eqn}) gives
\[(\phi')^2 \spd{y_1}{X} + 2(1-x^2) y_1 + 2 y_0 y_1 = - 2 \phi'
\mpd{y_0}{x}{X} - \phi'' \pd{y_0}{X}.\]
The homogeneous version of this equation is satisfied by both 
\[ \pd{y_0}{A} \quad \mbox{ and } \quad \pd{y_0}{\mu}.\]
Thus, by the Fredholm alternative, in order for there to be a
solution for $y_1$ we have the solvability conditions
\beqa
 2 \phi'
\int_0^1\mpd{y_0}{x}{X} \pd{y_0}{A}\, \d X  + \phi'' \int_0^1
\pd{y_0}{X} \pd{y_0}{A} \, \d X & = & 0,\\
 2 \phi'
\int_0^1\mpd{y_0}{x}{X} \pd{y_0}{\mu}\, \d X  + \phi'' \int_0^1
\pd{y_0}{X} \pd{y_0}{\mu} \, \d X & = & 0.\label{fred2}
\eeqa
These equations are our differential equations for $A$ and $\mu$ as
functions of $x$.

To enable us to write down these differential equations explicitly,
let us define  the function $Y(X,x,\Phi,A)$ by
\beqa
 \Phi \int_{Y_1(A,x)}^{Y}\frac{\d y}{c(A,x,y)^{1/2}}  &=&
 X  \label{2a}, \quad 0<X<X^*(A,x,\Phi) \\
X^*(A,x,\Phi) + \Phi \int^{Y_2(A,x)}_{Y}\frac{\d y}{c(A,x,y)^{1/2}}
&=& X  \label{2b}, \quad X^*(A,x,\Phi)<X<2X^*(A,x,\Phi), \qquad
\eeqa
and the function $\Phi(A,x)$ by
\beq
 \Phi(A,x) = \left(2\int_{Y_1(A,x)}^{Y_2(A,x)} \frac{\d y}{c(A,x,y)^{1/2}}
 \right)^{-1},\label{3} 
\eeq
where
\beqa
 c(A,x,y) &=&  A + 2 y - 2(1-x^2)y^2
- \frac{2}{3}y^3,\label{c}\\
X^*(A,x,\Phi)& = & \Phi \int_{Y_1(A,x)}^{Y_2(A,x)}\frac{\d
   y}{c(A,x,y)^{1/2}}.
\eeqa
Then $y_0 = Y(X+\mu(x),x,\Phi(A(x),x),A(x))$. Note that (\ref{3}) gives
 $X^*(A,x,\Phi(x,A))=1/2$ as we would expect, since
 $\Phi(A(x),x)=\phi'$ was chosen 
 to make the period 1 in $X$. We now have 
\beqas
 \mpd{y_0}{x}{X} &=&  \mpd{Y}{x}{X} +
\mpd{Y}{\Phi}{X}\left(\pd{\Phi}{x} + A'\pd{\Phi}{A}\right) +
\mpd{Y}{A}{X}A' +   \spd{Y}{X}\mu',\\
\pd{y_0}{A} & =& \pd{Y}{A} + \pd{Y}{\Phi}\pd{\Phi}{A},\\
\pd{y_0}{\mu} & =& \pd{Y}{X},\\
\phi'' & = &\pd{\Phi}{x} + A'\pd{\Phi}{A},
\eeqas
and the solvability conditions are
\beqa
\lefteqn{ 2 \Phi
\int_0^1\left( \mpd{Y}{x}{X} +
\mpd{Y}{\Phi}{X}\left(\pd{\Phi}{x} + A'\pd{\Phi}{A}\right) +
\mpd{Y}{A}{X}A' +   \spd{Y}{X}\mu'\right) \left(\pd{Y}{A} +
\pd{Y}{\Phi}\pd{\Phi}{A}\right)\, \d X  }\non \hspace{10cm}&&\\*
\mbox{ } + 
\left(\pd{\Phi}{x} + A'\pd{\Phi}{A}\right)\int_0^1 \pd{Y}{X}
\left(\pd{Y}{A} + \pd{Y}{\Phi}\pd{\Phi}{A}\right) \, \d X & = &
0,\hspace{2cm}\label{mu}\\  
\lefteqn{ 2 \Phi
\int_0^1\left( \mpd{Y}{x}{X} +
\mpd{Y}{\Phi}{X}\left(\pd{\Phi}{x} + A'\pd{\Phi}{A}\right)+
\mpd{Y}{A}{X}A' +   \spd{Y}{X}\mu'\right)
\pd{Y}{X}\, \d X}\non \hspace{10cm}&&\\
\mbox{ }  + \left(\pd{\Phi}{x} + A'\pd{\Phi}{A}\right) \int_0^1 
\left(\pd{Y}{X}\right)^2 \, \d X & = & 0.\hspace{2cm}\label{A}
\eeqa
Because of periodicity, any terms which are exact derivatives in $X$
will integrate to zero. 
Moreover, $Y$ and its derivatives with respect to $x$, $\Phi$ and $A$ are even,
while $Y_X$ and derivatives with respect to $x$, $\Phi$ and $A$ are
odd. 
Thus in 
fact equation (\ref{mu}) reduces to
\beq
\mu'=0,\label{mu3}
\eeq
so that $\mu$ is in fact constant. Since $Y$ is periodic in $X$ with
period $1$ we may take $\mu \in [0,1)$ without loss of generality. 
Equation (\ref{A}) can be written as
\beqa
\Phi \pd{\ov{Y_X^2}}{x} + \Phi \left(\pd{\Phi}{x} +
A'\pd{\Phi}{A}\right)  \pd{\ov{Y_X^2}}{\Phi} + \Phi A' \pd{\ov{Y_X^2}}{A}
+  \left(\pd{\Phi}{x} +
A'\pd{\Phi}{A}\right) \ov{Y_X^2} &=& 0,\label{A1}
\eeqa
where
\[ \ov{Y_X^2} = \int_0^1 \left(\pd{Y}{X}\right)^2 \, \d X.\]
Remarkably, eqn (\ref{A1}) is simply
\beq
\Phi \fdd{ \ov{Y_X^2}}{x}+  \fdd{\Phi}{x}\ov{Y_X^2} =
\fdd{}{x}\left(\Phi \ov{Y_X^2}\right)= 0,\label{Asol1}
\eeq
so that 
\beq
\Phi \ov{Y_X^2} = \mbox{constant} = k,\label{Asol}
\eeq
say.
Now, since
\[
\Phi^2 \left(\pd{Y}{X}\right)^2   =  
 A+2Y - 2(1-x^2)Y^2-\frac{2}{3}Y^3,
\]
\beqa
\Phi \ov{Y_X^2} &=&\frac{1}{\Phi} \int_0^1 \left( A + 2 Y - 2(1-x^2)Y^2
- \frac{2}{3}Y^3\right)\, \d X \non \\
& = & 2\int_{Y_1(A,x)}^{Y_2(A,x)} \left( A + 2 y - 2(1-x^2)y^2
- \frac{2}{3}y^3\right)^{1/2}\, \d y = k.\label{Asolution}
\eeqa
This is an implicit solution for $A(x)$, and is equivalent to equation
(54) in \cite{wong2008}, which was determined  by other means.
In fact, it is nothing more than the principle of adiabatic invariance
\cite{LandauLifshitz}, which states that a trajectory can be
approximated by moving slowly from one closed orbit to another as $x$
varies (closed
in terms of the fast scale $X$ when treating $x$ as constant) in such a 
way that the enclosed area remains constant.
We see that equation (\ref{Asol1}) follows directly from (\ref{fred2}) once we
realise that $\pd{y_0}{\mu} = \pd{y_0}{X}$.

\begin{figure}
\centering
\includegraphics[height=8cm]{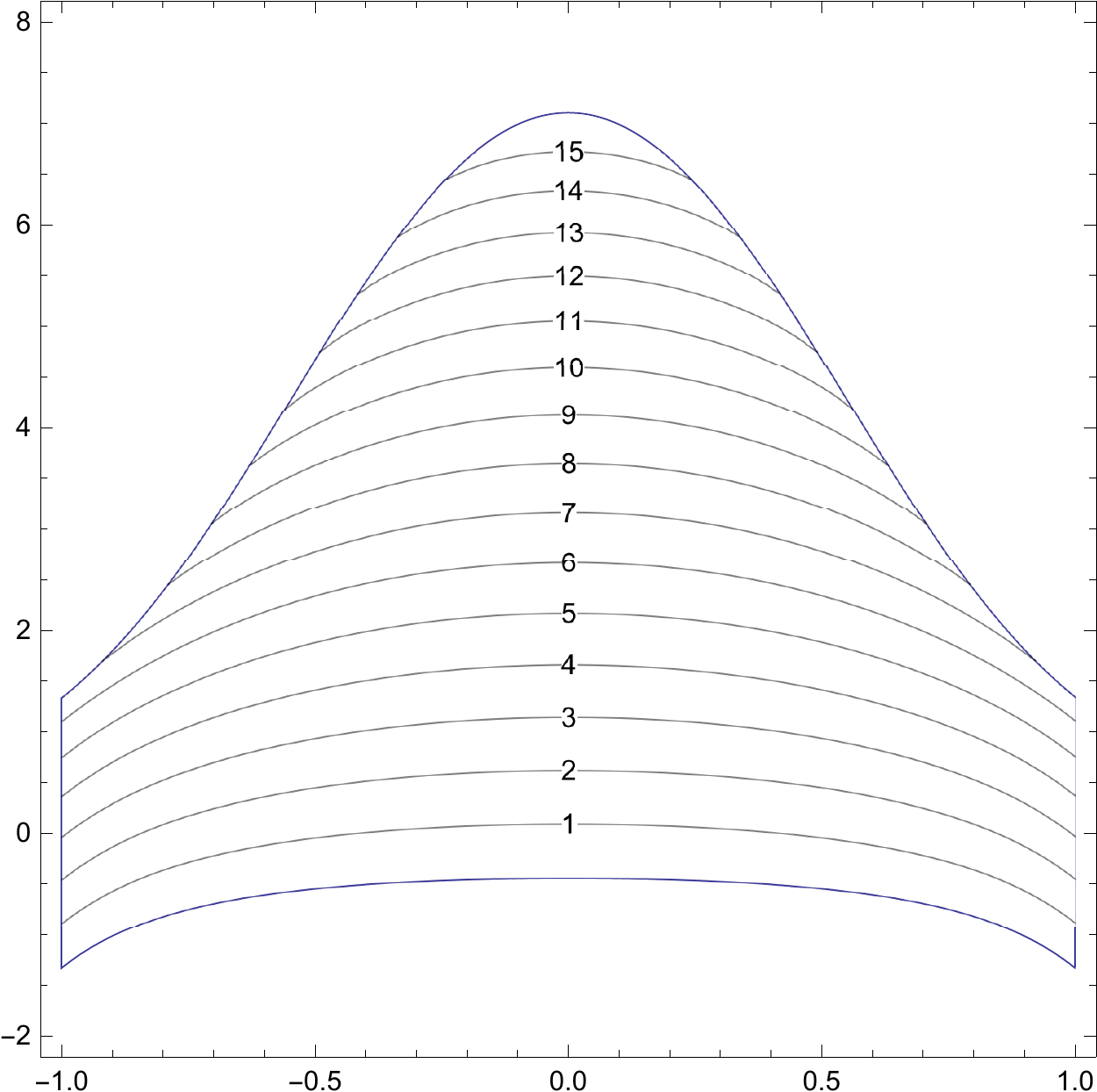}
\caption{$A$ as a function of $x$ for different values of $k$. The
  upper and lower curves are $A_1(x)$ and $A_2(x)$ respectively. For
  $0<k<k_1 \approx 6.78823$ the solution $A(x)$ stays between $A_1$
  and $A_2$ for all $-1<x<1$. For $k>k_1$ there are  turning points at
 $x=\pm x^*$ at which  
  $A(\pm x^*) = A_1(\pm x^*)$.}
\label{fig:Asoln}
\end{figure}

In Fig.~\ref{fig:Asoln} we show $A$ as a function of $x$ for integer
values of $k$ between $1$ and $15$. For
  $0<k<k_1 \approx 6.78823$ the solution (\ref{Asolution}) 
is valid in the whole domain $-1<x<1$. However, we will see in
\S\ref{sec:turningpoints} that 
 for $k>k_1$ there are turning points at  $x=\pm x^*$ at which  
  $\Phi=0$. For $|x|>x^*$ the solution will cease to be oscillatory
  and will instead be described by the outer solution
  (\ref{outersol}). Thus values $0<k<k_1$ correspond to solutions in
  which the spikes fill the domain, while values $k>k_1$ correspond to
  solutions in which there is an interior region near $x=0$ containing
  spikes, separated from the boundary layers by the  spike-less outer solution.

\begin{figure}
\centering
\includegraphics[height=8cm]{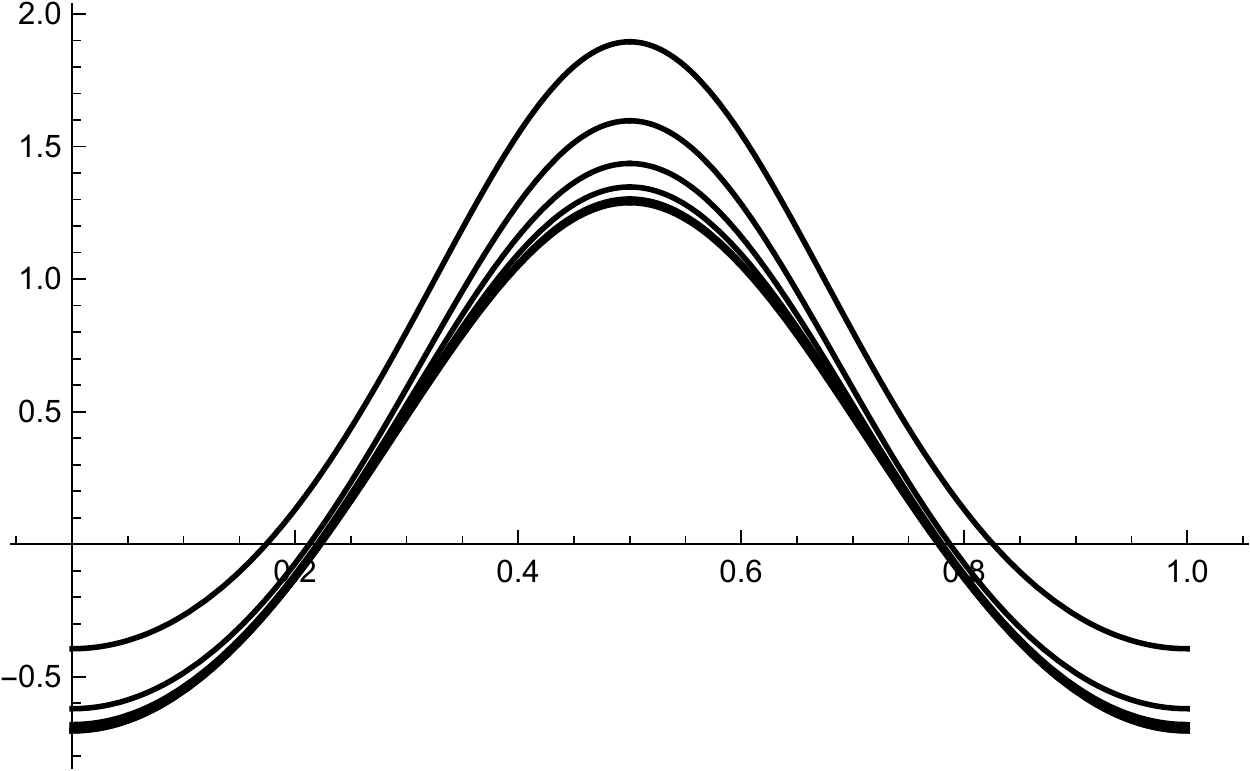}
\caption{The function $Y(X,x,\Phi(x,A(x)),A(x)$  for $k=5$ at $x=0$
  (bottom), $0.2$, $0.4$, $0.6$, $0.8$, and $1$ (top). Because of the
  leading-order nonlinearity in (\ref{maineqn}), the underlying oscillator
  (\ref{leading}) is nonlinear, so that the form of the oscillation (and not
  just the amplitude) varies with position.} 
\label{fig:Yexample}
\end{figure}

In Fig.~\ref{fig:Yexample} we show 
 $Y$ as a function of $X$ when $k=5$ and $A$ is given by 
(\ref{Asolution}), for $x = 0$, $0.2$, $0.4$, $0.6$, $0.8$, $1$. This
illustrates the fact that the form of the oscillation varies with
position, which is due to the leading-order nonlinearity in
(\ref{maineqn}), which causes the underlying oscillator
(\ref{leading}) to be nonlinear.

We have now determined the leading-order solution, up to the imposition
of the boundary conditions. To summarise, we have 
\beq
y_0 = Y\left(\frac{\phi(x)}{\eps}+\mu,x,\Phi(x,A(x)),A(x)\right),
\qquad \phi'(x)  =  \Phi(A(x),x), 
\eeq
where $Y(X,x,\Phi,A)$ is the function given by (\ref{2a})-(\ref{2b}) and 
\beqas
 \Phi & =& \left(2\int_{Y_1(A,x)}^{Y_2(A,x)}
  \frac{\d y}{(A + 
    2y - 2(1-x^2)y^2 - 2 
  y^3/3)^{1/2}} \right)^{-1},\\
k & = & 2\int_{Y_1(A,x)}^{Y_2(A,x)} \left( A + 2 y - 2(1-x^2)y^2
- \frac{2}{3}y^3\right)^{1/2}\, \d y,
\eeqas
where $\mu$ and $k$ are constants to be determined by the boundary
conditions. Note that we can choose $\phi(0)=0$ without loss of generality.

\subsection{Turning Points}
\label{sec:turningpoints}
Before we investigate the imposition of the boundary conditions to
determine the remaining unknown constants, we first discuss in more
detail the turning points that may appear in the solution.
These occur whenever \mbox{$\phi'(x) =0$}, i.e.
$\Phi(A(x),x) = 0$. Looking at (\ref{3}), we see this will happen when 
$Y_0 \ra Y_1$, that is, the left-most roots of $c(y)$ coalesce. (We
might also expect something strange to happen when the right-most
roots of $c(y)$ coalesce, that is, when $Y_1 \ra Y_2$. In that case,
however, because the range of integration in (\ref{3}) shrinks to zero
at the same rate at which the integrand blows up, $\Phi$ remains
finite. The case $Y_1 \ra Y_2$ corresponds to the limit in which the
constant $k \ra 0$.)

Thus, at the turning point, we have a double root of $c$, so that
there is a simultaneous root of $c$ and $\pd{c}{y}$. 
Such a double root occurs when $A = A_1(x)$ or $A=A_2(x)$ where
\beqas
A_1(x) & = & \frac{2}{3}\left(5 - 9 x^2 + 6 x^4 - 2 x^6 + 2(2  - 
    2 x^2  +  x^4 )^{3/2}\right),\\
A_2(x) & = & \frac{2}{3}\left(5 - 9 x^2 + 6 x^4 - 2 x^6 - 2(2 - 
   2 x^2  +  x^4  )^{3/2}\right).
\eeqas
When $A(x) = A_1(x)$ the left most roots of $c(y)$ coalesce ($Y_0 \ra
Y_1$). At $A(x) = A_2(x)$ the right most roots of $c(y)$ coalesce
($Y_1 \ra Y_2$). 
The functions
$A_1(x)$ and $A_2(x)$ form the upper and lower boundaries in
Fig. \ref{fig:Asoln} respectively.

The limiting value of $k$ for which there are no turning points in
$[-1,1]$ and the oscillations are present all the way to the boundary
is given by $A(1) = A_1(1)$, that is, it is the value of $k$ for which the
turning points lie at $x = \pm 1$.
We find this is 
\[ k = k_1 \approx 6.78823 .\]
 For  $k>k_1$ there will be an interior region between the two turning
 points $x=\pm x^*$ in which the solution is rapidly oscillating
 (i.e.~in which there are spikes), separated from the boundary layers
 by the  spike-less outer solution (\ref{outersol}).

\subsection{Boundary conditions}

Let us now seek to determine the remaining unknown constants $k$ and
$\mu$ by imposing the boundary conditions on our asymptotic solution.
We first look for solutions in which there are no turning points,
that is, in which the Kuzmak
approximation we have derived is valid all the way to the boundary.
The conditions $y(-1) = y(1) = 0$ imply
\beq
 Y\left(\frac{\phi(\pm 1)}{\eps}+\mu,\pm 1,\Phi(\pm 1,A(\pm 1)),A(\pm
1)\right) = 0.\label{bc}
\eeq
Now, the function $Y(X,\pm 1,\Phi(\pm 1,A(\pm 1)),A(\pm
1))$ has two zeros in the unit cell $0<X<1$ (see
Fig.~\ref{fig:Yexample}). Let us denote the smaller by $ 
X_0$; the larger is then given by $1-X_0$.  
Then,  at leading order, (\ref{bc}) gives
\beq
 \frac{\phi(1)}{\eps} + \mu = n\pm  X_0, \qquad \frac{\phi(-1)}{\eps}
+ \mu = m \pm  X_0,\label{bc1}
\eeq
where $n$, $m \in \Z$. 
These are two equations for the two unknown constants $k$ and
$\mu$. 

Eliminating $\phi(1)$,  noting that $\phi(1) = -\phi(-1)$, gives the
four possibilities
\beq
  n+m  = \quad 2\mu, \quad 2 \mu, \quad -2X_0 + 2 \mu, \quad 2X_0 + 2
  \mu,\label{nm}
\eeq
corresponding to choosing the signs in (\ref{bc1}) as $+-$, $-+$, $++$
and $--$ respectively.
In the first two cases (which give a symmetric solution) we must have
$\mu=0$ (corresponding to a solution with a minimum at the origin) or
$\mu=1/2$ (corresponding to a solution with a maximum at the origin).
We first analyse these symmetric solutions, and their associated
bifurcations, before returning to consider the non-symmetric solutions.

\subsection{Symmetric Solutions}
We consider here  solutions in which $\mu=0$   or
$\mu=1/2$. In this case equations (\ref{bc1}) reduce to 
\beq
\frac{\phi(1)}{\eps}  = n \pm  X_0
\qquad \mbox{ and } \qquad
\frac{\phi(1)}{\eps}  = n+ \frac{1}{2} \pm  X_0,\label{bc1a}
\eeq
respectively. We need to find the values of $k$ for which one of these
equations is satisfied.

\begin{figure}
\centering
\includegraphics[height=8cm]{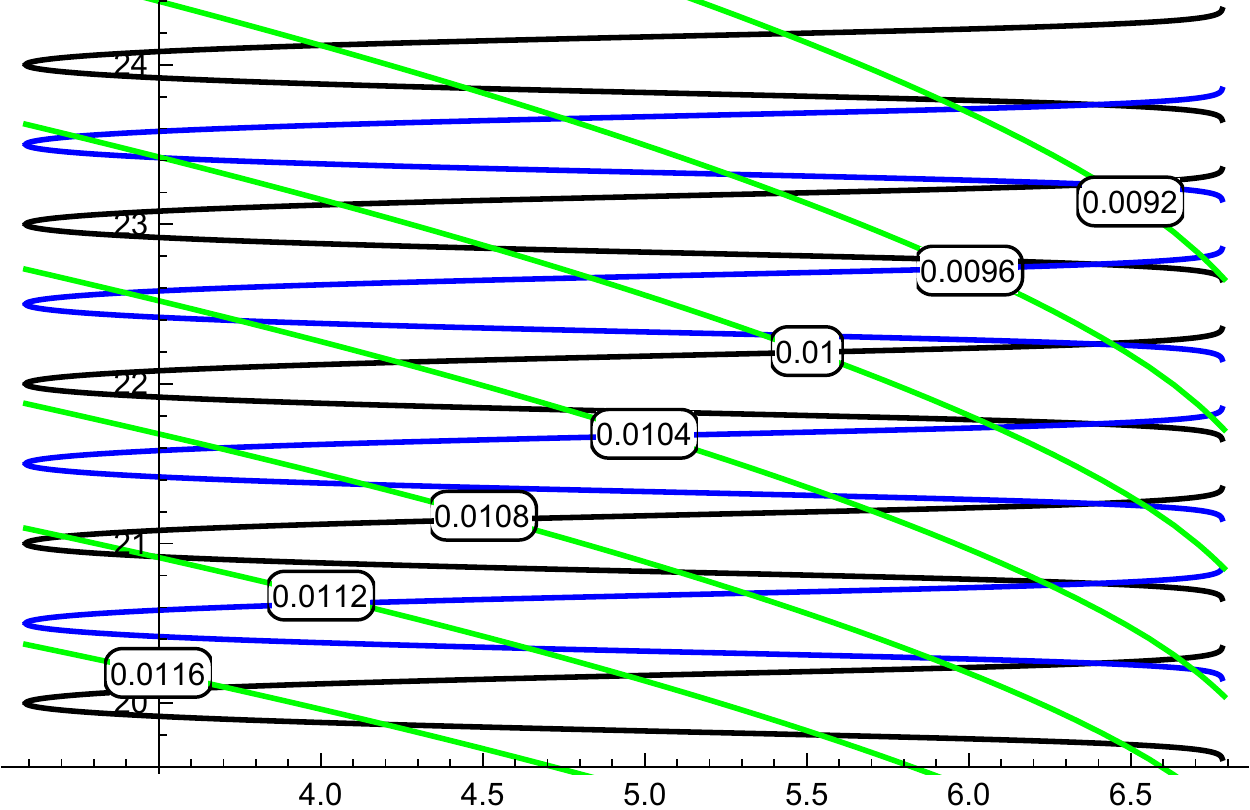}
\caption{The curves $n \pm X_0$ (black) and $n+1/2\pm X_0$ 
(blue) as a function of $k$ for $n$ in the range $20$ to $24$.
Also shown are the curves $\phi(1)/\eps$ for values of $\eps$
indicated (green). The intersections between
these curves give symmetric solutions of the problem.}
\label{fig:intersect}
\end{figure}

We show in Figure \ref{fig:intersect} the curves
\[ n \pm X_0 \quad \mbox{ (black)}\qquad \mbox{and} \qquad n
+\frac{1}{2}\pm X_0 
\quad\mbox{ (blue)}\]
as a function of $k$ for $n$ in the range $20$ to $24$.
The noses of these curves correspond to the value of $k$ at which $Y=
Y_X=0$ at $x = \pm 1$, which means that $A(\pm 1) = 0$.
This corresponds to
\[ k = k_0 = \frac{16\times 3^{3/4} \pi^{3/2}}{5 \Gamma(1/4)^2} \approx 3.08997 .\]
Also shown in Figure  \ref{fig:intersect} (green) are the curves
$\phi(1)/\eps$ for $\eps = 0.01 + 0.004 j$ 
with $j$ ranging from $-2$ to $4$. For a given value of $\eps$, the
intersections between the corresponding green curve and the black and
blue  curves give solutions of (\ref{bc1a}), and therefore correspond to
symmetric solutions of (\ref{maineqn}).  
We see that as $\eps$ decreases there is a succession of fold
bifurcations as a new pair of intersection points appears near $k =
k_0$. 

Initially, as $\eps$ decreases, one of the two new intersection points
moves to the left and one to the right. However, the left-moving point
soon reaches the nose of the blue/black curve, after which both
intersection points move to the right, in the direction of increasing $k$.
Eventually both approach  the limiting value $k =
k_1$, at which point a turning point appears near the boundary, and
another analysis takes over, since the boundary conditions should not
then be applied to the Kuzmak solution.  
After the turning point appears these solutions transition into
solutions whose oscillations do not
encompass the whole domain, but are restricted some smaller interval.

In Figure \ref{fig:bif} we show the figure analogous to Fig.~\ref{fig:intersect}
for  $\eps = 0.0335$. We illustrate the four intersection points,
along with the corresponding asymptotic approximation of the solutions.

\subsubsection{Maximum number of spikes}
The number of maxima in each solution is the number of complete periods
in $X$, which is  
\[ \frac{\phi(1)- \phi(-1)}{\eps} =
\frac{2\phi(1) }{\eps}.\]
Since $\phi(1)$ is monotonically decreasing in $k$ (see
Fig.~\ref{fig:intersect}), the largest value
of $\phi(1)$ occurs for $k = k_0$. This gives
\[ 2\phi(1) \approx 0.472537.\]
Thus the maximum number of spikes is
\[ \left\lfloor\frac{0.472537}{\eps}\right\rfloor \]
where $\lfloor x \rfloor$ denotes the largest integer less than $x$,
in agreement with the results of \cite{wong2008}.

\subsubsection{Proportion of solutions with oscillations filling the
  domain} 
The description in the introduction indicated that we might gradually add
spikes into the 
interior of the domain until there is no room to fit any more. 
Thus we might have expected that the proportion of solutions
containing turning points (i.e.~the proportion containing some
spike-free region) should tend to 1 as $\eps \ra 0$. However, minimum
number of spikes which may be present in a solution which does not
have turning points is given by the minimum value of 
\[ \left\lfloor\frac{2\phi(1) }{\eps}\right\rfloor\]
for $k_0<k<k_1$, which occurs when $k=k_1$, and is
\[ \left\lfloor\frac{0.415}{\eps}\right\rfloor.\]
Only for solutions with fewer spikes will the oscillations not fill
the domain.
Thus the proportion of solutions which do not have turning points is
approximately $0.12$
in the limit as $\eps \ra 0$.

\subsubsection{Position of the fold bifurcations}

The fold bifurcation is not exactly at  the nose of the blue or black
curve in Fig. \ref{fig:intersect}, though it approaches it as
$\eps \ra 0$.
At the bifurcation point the green and blue/black curves are tangent,
so that 
\beq \frac{1}{\eps}\fdd{\phi(1)}{k} = -\fdd{X_0}{k},\label{fold1}
\eeq
which must be satisfied at the same time as
\beq \frac{2\phi(1)}{\eps}  = n-  2X_0\label{fold2}
\eeq
(minus sign because the tangency is on the lower branch of the
blue/black curve).
Equations (\ref{fold1}) and (\ref{fold2}) form  two equations for $k$ and $\eps$
as a function of $n$. 
\begin{figure}
\centering
\includegraphics[height=8cm]{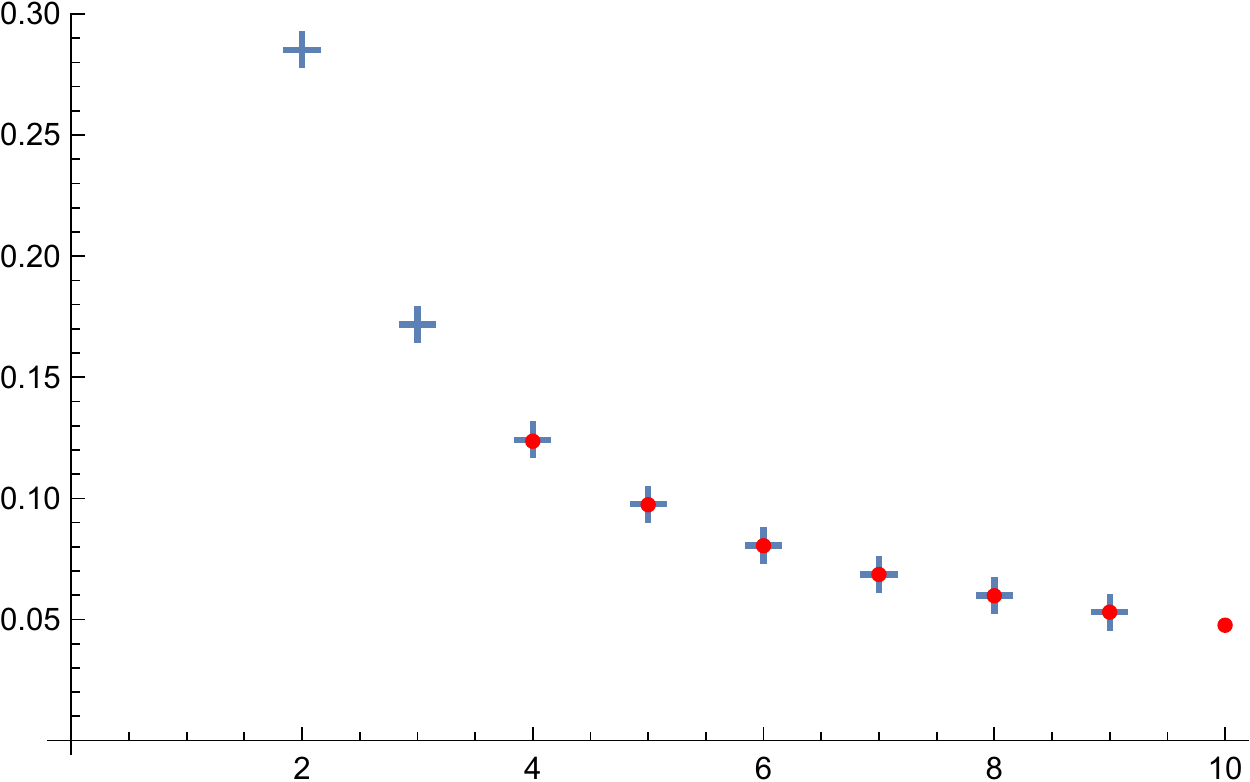}
\caption{The values of $\eps$ at the fold bifurcations, as a function
  of $n$. Equations (\ref{fold1}) and (\ref{fold2}) are valid in the limit
  $\eps \ra 0$, corresponding to $n \ra \infty$. The first value of
  $n$ for which a tangency point exists is $n=4$, which
  corresponds to the third fold bifurcation. Crosses indicate
  numerical results, circles the asymptotic approximation.} 
\label{fig:folds}
\end{figure}

We can approximate (\ref{fold1}), (\ref{fold2}) in the limit of small
$\eps$ 
(large $n$),
taking advantage of the fact that the bifurcation point is close to $k
=k_0$. Setting $k = k_0 + \delta$ we have
\[ X_0^2 \sim a \delta+\cdots \qquad \mbox{ where }\qquad 
 a =  \left.\fdd{X_0^2}{k}\right|_{k=k_0}.
\]
Equation (\ref{fold1}) gives
\[
- \frac{a^{1/2}}{2 \delta^{1/2}} \sim
\frac{1}{\eps}\left.\fdd{\phi(1)}{k}\right|_{k=k_0}, \] 
so that
\[ \delta \sim \frac{a
  \eps^2}{4}\left(\left.\fdd{\phi(1)}{k}\right|_{k=k_0}\right)^{-2}\]
Since
\[
\phi(1) \sim  \left.\phi(1)\right|_{k=k_0} + \delta
\left.\fdd{\phi(1)}{k}\right|_{k=k_0}+\cdots.\] 
equation (\ref{fold2}) now gives
\[
\frac{2 \left.\phi(1)\right|_{k=k_0}}{\eps} + 
 \frac{a
  \eps}{2}\left(\left.\fdd{\phi(1)}{k}\right|_{k=k_0}\right)^{-1} \sim n
+  a  \eps\left(\left.\fdd{\phi(1)}{k}\right|_{k=k_0}\right)^{-1}
\]
Thus
\beqa
\eps &\sim& \frac{2 \left.\phi(1)\right|_{k=k_0}}{ n
+  \frac{a
  \eps}{2}\left(\left.\fdd{\phi(1)}{k}\right|_{k=k_0}\right)^{-1}}
\sim\frac{2 \left.\phi(1)\right|_{k=k_0}}{ n
+  
 \frac{a
   \left.\phi(1)\right|_{k=k_0}}{n}
\left(\left.\fdd{\phi(1)}{k}\right|_{k=k_0}\right)^{-1}}\non \\
&
\approx&\frac{0.472537}{ n
-
 \frac{0.8344}{n}}.\label{asyfold}
\eeqa

In Figure \ref{fig:folds} we show  the values of $\eps$ at the fold
bifurcations as a function
  of $n$.  The first value of
  $n$ for which a tangency point exists is $n=4$, which
  corresponds to the third fold bifurcation. Although the asymptotic
  approximation is valid as $n \ra \infty$, the agreement with the
  numerical solution is
  remarkably good even at small $n$.

\subsection{Non-symmetric solutions}
Let us now return to consider the other solutions of (\ref{bc1}). 
From (\ref{nm}) we see that for non-symmetric solutions we must have
 $\mu = \pm X_0$ (modulo 1) or  $\mu = 1/2 \pm X_0$ (modulo 1). 
In each case we find $n=-m$. For $\mu = \pm X_0$ we find  that
(\ref{bc1}) becomes
\beq \frac{\phi(1)}{\eps} = n.\label{one}
\eeq
Similarly, for $\mu = 1/2\pm X_0$ we find that (\ref{bc1}) becomes
\beq \frac{\phi(1)}{\eps} = n-1/2.\label{half}
\eeq
Using Fig.~\ref{fig:intersect} to illustrate these solutions, 
we see that they correspond to the intersection points between the green
curves $\phi(1)/\eps$ and the horizontal
lines $n$ and $n-1/2$. The value of $\mu$ then corresponds to the
vertical distance from this intersection point to the nearest black
curve (modulo 1). Each intersection point gives two asymmetric
solutions, corresponding to the two values $\mu =  \pm X_0$ in case
(\ref{one}) or $\mu = 1/2
\pm X_0$ in case (\ref{half}).

In Figure \ref{fig:bifa} we show these intersection points and
corresponding $\mu$ values for  $\eps = 0.0335$. 
There are two intersection points, each corresponding to two
solutions. 
The corresponding asymptotic approximation of the solutions is also
illustrated.

When $X_0 = 0$ we find $\mu = 0$ or $\mu = 1/2$, corresponding to a
pitchfork bifurcation  at which the non-symmetric solutions bifurcate
from one of the symmetric 
branches.

\subsubsection{Position of the pitchfork bifurcations}

The pitchfork bifurcation occurs when the intersection point between
the green and blue/black curves lies exactly at the nose of those
curves. At that point 
\beq X_0 = 0 \qquad \mbox{ and } \qquad \frac{2\phi(1)}{\eps}  =n.
\label{pitch}  
\eeq 
The first equation gives $k = k_0 \approx 3.08997$. Then $\phi(1)$ is
determined and the second equation gives
\beq
 \eps = \eps_n = \frac{2\phi(1)}{n}\approx \frac{0.472537}{n}.\label{asypitch}
\eeq

\begin{figure}
\centering
\includegraphics[height=8cm]{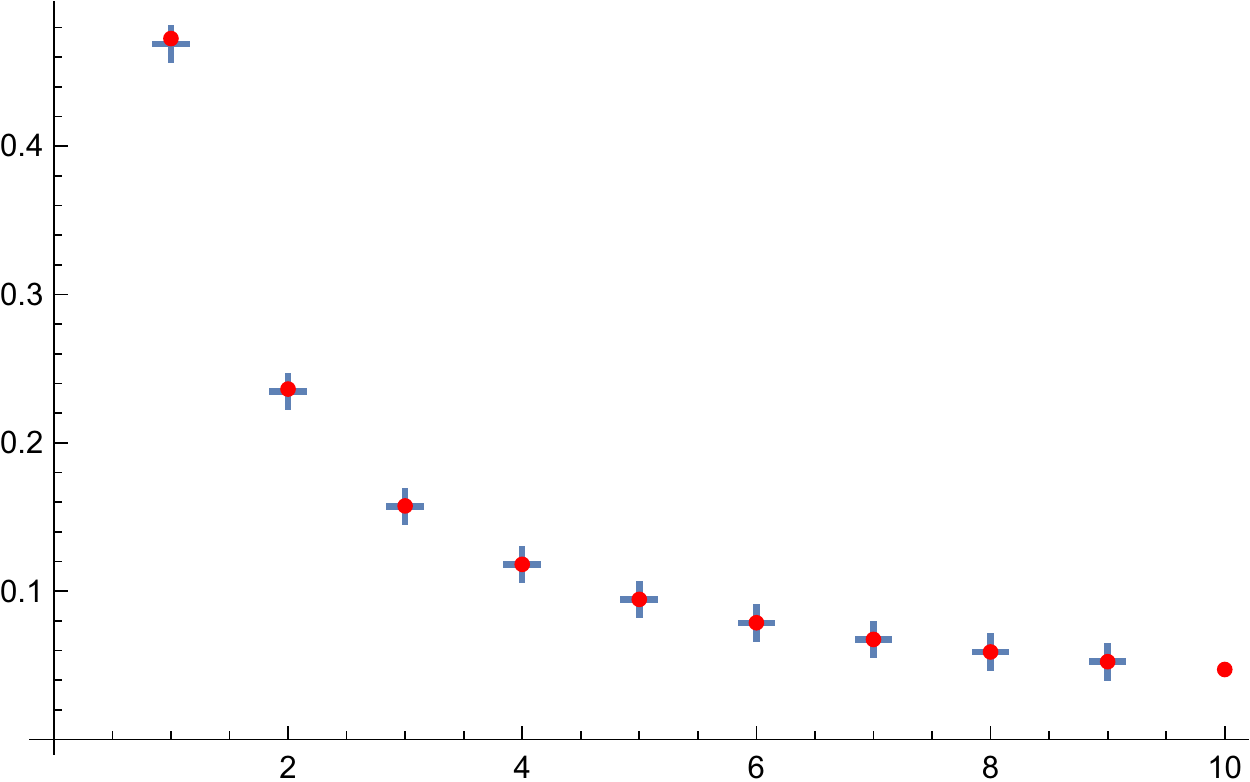}
\caption{The values of $\eps$ at the pitchfork bifurcations, as a function
  of $n$. Equation (\ref{pitch}) is valid in the limit
  $\eps \ra 0$, corresponding to $n \ra \infty$.  Crosses indicate
  numerical results, circles the asymptotic approximation.} 
\label{fig:pitch}
\end{figure}

In Figure \ref{fig:pitch} we show  the values of $\eps$ at the pitchfork
bifurcations as a function
  of $n$.  Although the asymptotic
  approximation is valid as $n \ra \infty$, the agreement with the
  numerical solution is
  remarkably good even at small $n$.

From (\ref{asyfold}) and (\ref{asypitch}) we see that the separation
between the fold and pitchfork bifurcations is approximately
\[ \frac{0.3943}{n^3}\]
as $n \ra \infty$, and that the proportion of values of $\eps$ for
which there are $4N+2$ solutions rather than $4N$ solutions therefore
shrinks as $3.737 \eps^2$ as $\eps \ra 0$.

\begin{figure}
\centering
\scalebox{0.8}{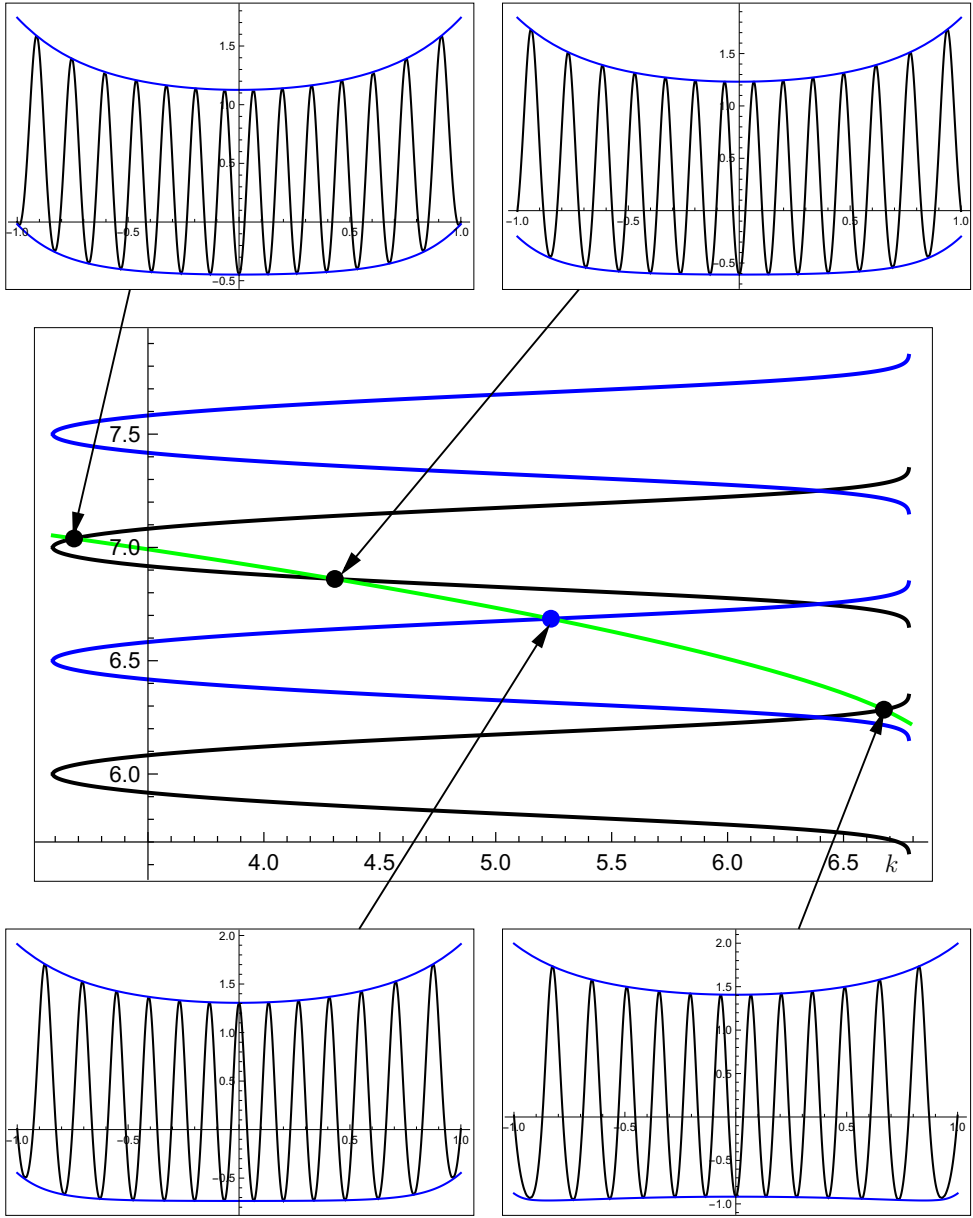}
\caption{{\bf Symmetric Solutions}. The curves $n \pm X_0$ (black) and $n+1/2\pm X_0$ 
(blue) as a function of $k$ for $n = 6$, $7$.
The green curve is $\phi(1)/\eps$ with $\eps = 0.0335$.
The intersection points, highlighted, correspond to symmetric
solutions of the problem. The lower and upper envelopes of the
individual solutions are the curves $Y_1(A(x),x)$ and $Y_2(A(x),x)$
respectively.} 
\label{fig:bif}
\end{figure}
\begin{figure}
\begin{center}
\scalebox{0.8}{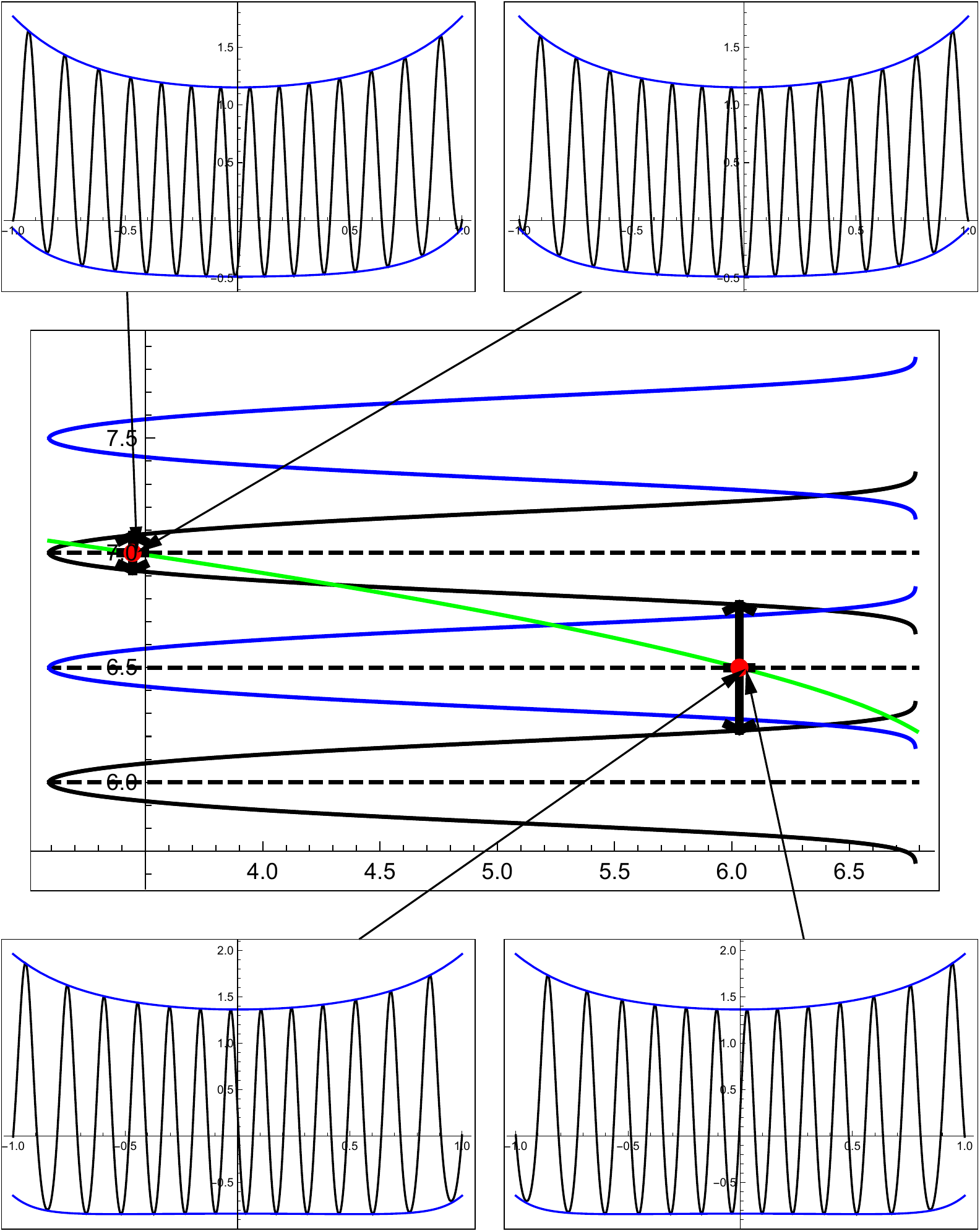}
\end{center}
\caption{{\bf Non-symmetric Solutions}. The curves $n \pm X_0$ (black)
  and $n+1/2\pm X_0$  
(blue) as a function of $k$ for $n = 6$, $7$.
The green curve is $\phi(1)/\eps$ with $\eps = 0.0335$.
The intersection points of the green curve with integer and
half-integer values, highlighted, correspond to non-symmetric
solutions of the problem. The distance of these intersection points
to the nearest black curve give the (signed) value of $\mu$. The lower
and upper envelopes of the 
individual solutions are the curves $Y_1(A(x),x)$ and $Y_2(A(x),x)$
respectively.} 
\label{fig:bifa}
\end{figure}

\subsection{Solutions with turning points}
Thus far we have  analysed solutions in which the Kuzmak approximation is valid
all the way up to the boundary, which enabled us to capture quite well
the structure of the bifurcation diagram shown in Figure
\ref{fig:bifurcation-diagram}. We now complete our asymptotic analysis
by considering those solutions with turning points, in which the
oscillations are confined to an interior region.

Suppose there are turning points at $\pm x^*$ (with $x^*>0$), so that $A(x^*) =
A_1(x^*)$. For $|x|>x^*$ the solution does not oscillate, and (outside
the boundary layers) $y_0$ is simply given by the outer solution
(\ref{outersol}):
\beq
 y_{\mathrm{out}}(x) = -1 + x^2 - \sqrt{x^4-2 x^2+2}, \qquad
 |x|>x^*.\label{outer} 
\eeq
In this case $k$ is determined  (at leading order) by the
condition that 
\beq
 \pd{Y}{X}=0\label{derivative}
\eeq
at $x = x^*$, with $Y$ negative, so that the oscillating solution can
join smoothly onto the 
non-oscillating solution. Note that 
continuity in the solution is automatic, since
at the turning point $Y = Y_1 = Y_0$, so that $Y$ is a root of both
the cubic $c(Y)$ and also its derivative.
 Since, from (\ref{2a})-(\ref{2b}),  $\pd{Y}{X}=0$ (and
$Y$ is negative) when $X=0$, equation (\ref{derivative}) gives 
\beq
 \frac{\phi(x^*)}{\eps}+\mu = n, \qquad  -\frac{\phi(x^*)}{\eps}+\mu =
 m.\label{jon}
\eeq
Thus we are forced to choose either $\mu=0$ or $\mu=1/2$: the central
part of the solutions is symmetric for all solutions with turning points.
Note that 
\[ \phi' \sim \frac{a}{\log (x-x^*)},\]
for some constant $a$
as the turning point is approached with the result that the separation
between spikes is $O(\eps \log  \eps)$ there, in agreement with the
analysis in \cite{macgillivray2000} on the two spike solution.

\begin{figure}
\begin{center}
\scalebox{0.7}{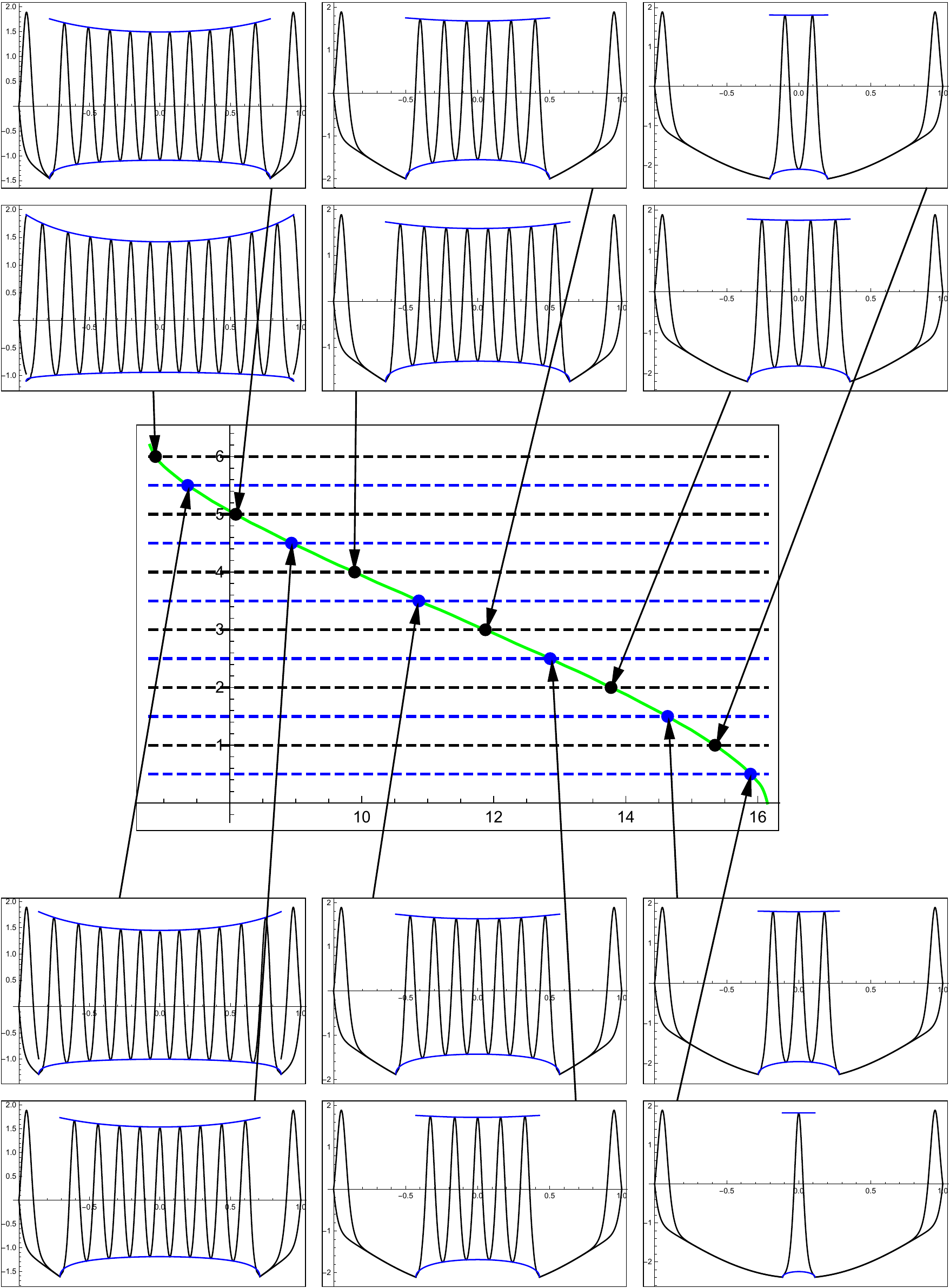}
\end{center}
\caption{{\bf Solutions with turning points}. 
The green curve is $\phi(x^*)/\eps$ with $\eps = 0.0335$ and $x^*$
satisfying $A(x^*) =A_1(x^*)$, shown as a function
of $k$.
The intersection points with integer and
half-integer values, highlighted, correspond to solutions of the
problem. Each point gives four distinct solutions, corresponding to
the four combinations of boundary layers at $x=\pm1$.  The lower
and upper envelopes of the 
individual solutions are the curves $Y_1(A(x),x)$ and $Y_2(A(x),x)$
respectively. Each solution plot shows both boundary layer
possibilities at each end. For $n=6$ and $n=5.5$ the turning point is
in the boundary layer, with the result that our approximation there is
inaccurate.} 
\label{fig:biftp}
\end{figure}

As described in the introduction, the outer solution (\ref{outer})
does not satisfy the boundary conditions 
at $x = \pm 1$, where there are boundary layers.
A uniform approximation valid for $|x|>x^*$
is 
\beqas
 y& \sim& 
-1 + x^2 - \sqrt{2 - 2 x^2 + x^4} + 
 3\, \sech^2\left( \pm \frac{(1 - x)}{\eps\sqrt{2}} +
   \tanh^{-1}\left(\sqrt{\frac{2}{3}}\right) \right)\\
&& \mbox { }+ 
 3\, \sech^2\left( \pm \frac{(1 + x)}{\eps\sqrt{2}}+
   \tanh^{-1}\left(\sqrt{\frac{2}{3}}\right) \right).
\eeqas

In Figure \ref{fig:biftp} we show the solutions for $\eps=0.0335$, for
which $\phi(x^*)/\eps$ may take any integer or half-integer value up
to $n=6$. The discontinuity in the gradient of the solution at $x=\pm
x^*$ in these plots is due to the fact that we imposed continuity of
the derivative 
only at leading order in $\eps$ there, using (\ref{derivative}). Full
continuity of the derivative implies
\[ \fdd{y_{\mathrm{out}}}{x} = \fdd{Y}{x} = \frac{\Phi}{\eps}
\pd{Y}{X} + \pd{Y}{x} + \pd{Y}{\Phi}\left(\pd{\Phi}{x} + \pd{\Phi}{A}
  \fdd{A}{x}\right) + \pd{Y}{A} \fdd{A}{x},
\]
which would introduce an $O(\eps)$ correction into equations (\ref{jon}).
Each solution curve in Figure \ref{fig:biftp} shows four distinct
solutions overlaid, corresponding to 
the four combinations of boundary layers at $x=\pm1$.

These 48 solutions, together with the 8 solutions shown previously in
Figures \ref{fig:bif} and \ref{fig:bifa},
make up the 56 solutions to the problem when $\eps=0.0335$.

\section{Conclusion}

The computational and asymptotic analysis 
  we have 
presented gives a 
novel and complete taxonomy of the solutions of Carrier's problem
(\ref{maineqn}). 

Using deflated continuation we found a rather striking bifurcation
diagram, containing an apparently infinite
number of mutually disconnected components. Each component (except for the first two)
contains one fold bifurcation, at which two solutions of
(\ref{maineqn})   appear, 
and one pitchfork 
bifurcation, at which a further two solutions of (\ref{maineqn})
appear. Solutions on the same connected component have the  same
number of interior maxima.  For values of $\eps$ which do not lie between the
fold and pitchfork bifurcations of a connected component the number of
solutions of (\ref{maineqn}) is a multiple of 4, as claimed by 
 Bender \& Orszag \cite{bender1999}. However, between the fold and
 pitchfork bifurcations the number of solutions is $4n+2$ with $n \in \Z$.

Our asymptotic analysis used Kuzmak's method to construct approximate
solutions of (\ref{maineqn}). Both the fold and pitchfork bifurcation
points were predicted accurately.
We found that the separation between these bifurcation points tends
quickly to zero as $\eps \ra 0$, so that the proportion of values of
$\eps$ for which there are  $4n+2$ solutions rather than $4n$
solutions tends to zero as $3.737
\eps^2$ as $\eps
\ra 0$. 

We gave an alternative derivation of the result of Wong and Zhao that
the  maximum number of internal maxima 
is asymptotically
\[ \left\lfloor\frac{0.472537}{\eps}\right\rfloor. \]
Moreover, we found that approximately 12\% of solutions of the problem
have oscillations which fill the domain. The remaining 88\% of
solutions have oscillations in an interior region near $x=0$
separated from boundary layers by a non-oscillating outer solution.

The methods we have used are in no way specific to
(\ref{maineqn}). Carrier's problem provides a nice example, but any
slowly-varying phase plane with closed orbits would be amenable to our
approach.

\bibliographystyle{siam}
\bibliography{literature}

\appendix
\section{Parameter values of the initial bifurcations}

In the interest of completeness we tabulate the values of
$\eps$ at which the first four pitchfork and fold
bifurcations occur. The solution and parameter value at
which a simple bifurcation occurs satisfy an augmented system
of integro-differential equations \cite{moore1980}:
\begin{equation} \label{moore}
F(y, v, \eps) = \begin{bmatrix}
\eps^2 y'' + 2(1-x^2) y + y^2 - 1 \\
\eps^2 v'' + 2(1-x^2) v + 2yv \\
\|v\|^2 - 1
\end{bmatrix} = 0,
\end{equation}
where $y$ is the solution at the bifurcation point, $v$ is
the eigenfunction in the nullspace of the Fr\'echet
derivative of the equation, $\eps$ is the value of the
parameter at the bifurcation, and $\|\cdot\|$ denotes
the $L^2([-1, 1])$ norm.

As we wish to compute the parameter values to high accuracy,
a spectral discretization was chosen to approximate the
solutions of \eqref{moore}. We thus employed the Chebfun
system of Trefethen and co-workers
\cite{driscoll2014,birkisson2012}. Solving \eqref{moore} can
be rather difficult, and the main art in its solution is the
construction of good initial guesses for $(y, v, \eps)$.
These were computed as follows.

For each bifurcation, an initial guess $(\tilde{y},
\tilde{\eps})$ for the solution and parameter was acquired
from the data produced by deflated continuation. The finite
element solution $\tilde{y}$ was evaluated at 200 Chebyshev
points of the second kind and its Chebyshev interpolant
$\mathcal{I}\tilde{y}$ was constructed with Chebfun.
Carrier's problem at $\eps = \tilde{\eps}$ was then solved
with this initial guess, yielding $\hat{y}$, to ensure that
the first equation of \eqref{moore} had small residual. The
differential operator was linearized at $(\hat{y},
\tilde{\eps})$ to compute its eigenfunction $\hat{v}$ with
eigenvalue closest to zero; this ensured that the second and
third equations had small residual. The triplet $(\hat{y},
\hat{v}, \tilde{\eps})$ was then supplied as initial guess
to the solver for \eqref{moore}. The fold bifurcations
typically converged in four or five Newton iterations, while
the pitchfork bifurcations typically converged in twenty to
thirty iterations. In all cases Chebfun's error estimate for
the solution of \eqref{moore} was less than $10^{-10}$; no
further accuracy was possible due to the use of double
precision arithmetic.


\begin{table}[h!]
\centering
\begin{tabular}{c|c|c|c}
\toprule
Connected  & Computed  &Asymptotic & Relative  \\
component & $\eps$ & estimate & error\\
\midrule
1   & 0.46886251 & 0.472537& 0.007837  \\
2   & 0.23472529 & 0.236269& 0.006574 \\
3   & 0.15703946 & 0.157512& 0.003012 \\
4   & 0.11798359 & 0.118134& 0.001278 \\
\bottomrule
\end{tabular}
\caption{Computed parameter values for the first four pitchfork bifurcations. The asymptotic estimates are those of \eqref{asypitch}.}
\end{table}


\begin{table}[h!]
\centering
\begin{tabular}{c|c|c|c}
\toprule
Connected  & Computed  &Asymptotic & Relative \\
component & $\eps$ &  estimate & error \\
\midrule
2   & 0.28522538 &  0.298545 & 0.0467 \\
3   & 0.17186970 &  0.173608 & 0.01011\\
4   & 0.12421206 &0.124634 & 0.003397\\
5   & 0.09762446 &  0.0977706 & 0.001497\\
\bottomrule
\end{tabular}
\caption{Computed parameter values for the first four fold bifurcations. The asymptotic estimates are those of \eqref{asyfold}.}
\end{table}
\section{Approximation for large $\eps$}
For completeness we give here an asymptotic approximation to the two
solutions which continue to exist when $\eps$ is large. 
Expanding $y$ in an inverse power series in $\eps$ as
\[ y \sim y_0 + \eps^{-2} y_1  + \cdots,\]
gives at leading order
\[ y_0'' = 0, \qquad y_0(-1) = y_0(1) = 0,\]
with solution $y_0 \equiv 0$, indicating that $y$ is not $O(1)$ but
must be rescaled in some way.
Expanding 
\[ y \sim \eps^{-2} y_0 + \eps^{-4} y_1  + \cdots,\]
gives at leading order
\[ y_0'' = 1, \qquad y_0(-1) = y_0(1) = 0,\]
with solution 
\[ y_0 = \frac{x^2-1}{2}.\]
This solution has no internal maximum, and is the continuation to
large $\eps$ of the solution in panel $1$ of Fig \ref{fig:bifurcation-diagram}.

The second solution is
found by expanding $y$ as 
\[ y \sim \eps^2 y_0 + y_1  + \cdots,\]
to give at leading order
\[ y_0'' + y_0^2 = 0, \qquad y_0(-1) = y_0(1) = 0,\]
with solution 
\[ 1+ x = \frac{\sqrt{3}}{\sqrt{2}}\int_0^y \frac{\d u}{(y_{\mathrm{max}}^3 - u^3)^{1/2}},\]
where $y_{\mathrm{max}}$, the value of $y$ at $x=0$, satisfies
\[ 1 = \frac{\sqrt{3}}{\sqrt{2}}\int_0^{y_{\mathrm{max}}} \frac{\d
  u}{(y_{\mathrm{max}}^3 - u^3)^{1/2}} = \sqrt{\frac{3 \pi}{2 y_{\mathrm{max}}}}
\frac{\Gamma(4/3)}{\Gamma(5/6)},\]
so that
\[ y_{\mathrm{max}} = \frac{3 \pi \Gamma(4/3)^2}{2 \Gamma(5/6)^2}.\]
This solution has one internal maximum, and is the continuation to
large $\eps$ of the solution in panel $2$ of Fig \ref{fig:bifurcation-diagram}.

\end{document}

%% file: figures/bif.pdf_t
\begin{picture}(0,0)%
\includegraphics{bif.pdf}%
\end{picture}%
\setlength{\unitlength}{3947sp}%
\begingroup\makeatletter\ifx\SetFigFont\undefined%
\gdef\SetFigFont#1#2#3#4#5{%
  \reset@font\fontsize{#1}{#2pt}%
  \fontfamily{#3}\fontseries{#4}\fontshape{#5}%
  \selectfont}%
\fi\endgroup%
\begin{picture}(7599,9549)(889,-8923)
\put(7801,-6211){\makebox(0,0)[lb]{\smash{{\SetFigFont{12}{14.4}{\rmdefault}{\mddefault}{\updefault}{\color[rgb]{0,0,0}$k$}%
}}}}
\end{picture}%

%% file: figures/bifa.pdf_t
\begin{picture}(0,0)%
\includegraphics{bifa.pdf}%
\end{picture}%
\setlength{\unitlength}{3947sp}%
\begingroup\makeatletter\ifx\SetFigFont\undefined%
\gdef\SetFigFont#1#2#3#4#5{%
  \reset@font\fontsize{#1}{#2pt}%
  \fontfamily{#3}\fontseries{#4}\fontshape{#5}%
  \selectfont}%
\fi\endgroup%
\begin{picture}(7599,9549)(889,-8923)
\put(7801,-6211){\makebox(0,0)[lb]{\smash{{\SetFigFont{12}{14.4}{\rmdefault}{\mddefault}{\updefault}{\color[rgb]{0,0,0}$k$}%
}}}}
\put(6751,-4411){\makebox(0,0)[lb]{\smash{{\SetFigFont{12}{14.4}{\rmdefault}{\mddefault}{\updefault}{\color[rgb]{0,0,0}$\mu$}%
}}}}
\put(6751,-4861){\makebox(0,0)[lb]{\smash{{\SetFigFont{12}{14.4}{\rmdefault}{\mddefault}{\updefault}{\color[rgb]{0,0,0}$\mu$}%
}}}}
\put(2196,-3796){\makebox(0,0)[lb]{\smash{{\SetFigFont{12}{14.4}{\rmdefault}{\mddefault}{\updefault}{\color[rgb]{0,0,0}$\mu$}%
}}}}
\put(2181,-3611){\makebox(0,0)[lb]{\smash{{\SetFigFont{12}{14.4}{\rmdefault}{\mddefault}{\updefault}{\color[rgb]{0,0,0}$\mu$}%
}}}}
\end{picture}%

%% file: figures/biftp.pdf_t
\begin{picture}(0,0)%
\includegraphics{biftp.pdf}%
\end{picture}%
\setlength{\unitlength}{3947sp}%
\begingroup\makeatletter\ifx\SetFigFont\undefined%
\gdef\SetFigFont#1#2#3#4#5{%
  \reset@font\fontsize{#1}{#2pt}%
  \fontfamily{#3}\fontseries{#4}\fontshape{#5}%
  \selectfont}%
\fi\endgroup%
\begin{picture}(8425,11424)(1039,-10573)
\put(7208,-6478){\makebox(0,0)[lb]{\smash{{\SetFigFont{12}{14.4}{\rmdefault}{\mddefault}{\updefault}{\color[rgb]{0,0,0}$k$}%
}}}}
\end{picture}%